\documentclass[twoside]{amsart}
\usepackage{latexsym}
\usepackage{amssymb,amsmath,amsopn,amsthm,amsfonts,amscd}
\usepackage[dvips]{graphicx}   

\newtheorem{thm}{Theorem}
\newtheorem{lem}{Lemma}
\theoremstyle{definition}
\newtheorem{defn}{Definition}

\newtheorem{conj}{Conjecture}

\newtheorem{cor}{Corollary}
\renewcommand{\Re}{\mathbb R}
\renewcommand{\S}{\mathbb S}

\renewcommand{\S}{\mathbb{S}}

\newcommand{\eqdef}{\stackrel{\rm{def}}{=}}

\def\bea{\begin{eqnarray}}
\def\eea{\end{eqnarray}}
\def\ben{\begin{equation}}
\def\een{\end{equation}}
\def\brr{\begin{array}}
\def\err{\end{array}}

\usepackage{chngcntr}
\counterwithout{figure}{section}
\usepackage[dvips]{graphicx}
\usepackage{pdfpages}
\usepackage{caption}
\usepackage{subcaption}

\begin{document}
\title[Genealogy of convex solids]{A genealogy of convex solids via local and global bifurcations of gradient vector fields}
\author[G. Domokos, Z. P.J. Holmes \and Z. L\'angi ]{G\'abor Domokos, Philip Holmes  \and Zsolt L\'angi}
\address{G\'abor Domokos, Dept. of Mechanics, Materials and Structures, Budapest University of Technology,
M\H uegyetem rakpart 1-3., Budapest, Hungary, 1111}
\email{domokos@iit.bme.hu}
\address{Zsolt L\'angi, Dept. of Geometry, Budapest University of Technology,
Egry J\'ozsef u. 1., Budapest, Hungary, 1111}
\email{zlangi@math.bme.hu}
\address{Philip Holmes, Program in Applied and Computational Mathematics and Dept. of Mechanical and Aerospace Engineering, Princeton University, Princeton NJ 08544}
\email{pholmes@math.princeton.edu}




\begin{abstract}

Three-dimensional convex bodies can be classified in terms of the
number and stability types of critical points on which they can
balance at rest on a horizontal plane. For typical bodies these are
nondegenerate maxima, minima, and saddle-points, the numbers of which
provide a primary classification. Secondary and tertiary
classifications use graphs to describe orbits connecting these
critical points in the gradient vector field associated with each
body. In previous work it was shown that these classifications are
complete in that no class is empty. Here we construct 1- and
2-parameter families of convex bodies connecting members of adjacent
primary and secondary classes and show that transitions between them
can be realized by codimension 1 saddle-node and saddle-saddle
(heteroclinic) bifurcations in the gradient vector fields. Our results
indicate that all combinatorially possible transitions can be realized
in physical shape evolution processes, e.g. by abrasion of sedimentary
particles.

\end{abstract}

\maketitle

\tableofcontents

\vfill \newpage

\section{Introduction}
\label{s1}

\subsection{Motivation and background}
\label{s1.1}

The evolution of shapes of abrading bodies, such as pebbles in river
beds and on beaches, has been studied for over 70 years
(e.g. \cite{Rayleigh1, Rayleigh2,Rayleigh3,Firey,Bloore}). Data from
NASA's Curiosity Rover on Mars \cite{Mars1,Mars2} has rekindled
interest in the subject. In addition to classical shape indices such
as axis ratios and roundness \cite{Zingg, Illenberger},  a recent
approach considers the evolution of the number of \emph{static
equilibrium points} $N(t)$ on the surface of an abrading body, i.e.,
points on which the body can balance at rest on a horizontal plane
\cite{VD1,DSSZV,Plos,D}. Unlike shape indices, which require length
measurements, the integer $N(t)$ can be counted in simple  experiments
\cite{DSSZV}. 
 
Abrasion occurs primarily on a body's convex hull, so to formulate a
precise  and relatively simple model we restrict our analysis to
convex bodies $K$ of uniform density, with surfaces described by
scalar Euclidean distance functions $r_K$ measured from the center of
mass $C_K$. For such bodies static equilibria are critical points of
$r_K$ at which the gradient $\nabla r_K=0$.
 
The surface $\partial K$ of a generic convex body $K$ can exhibit
three types of nondegenerate critical points: local minima, maxima and
saddle-points, which are sinks, sources and saddles of the gradient
vector field $\mathbf{v} = \nabla r_K$. Let $S,U,H$ respectively denote
the number of each of these points. Since $\partial K$ is a
topological 2-sphere, the Poincar\'e-Hopf Theorem \cite{Arnold1}
implies that
\ben \label{Poincare}
S + U - H = 2 . 
\een

\begin{figure}[here]
\begin{center}
\includegraphics[width=120mm]{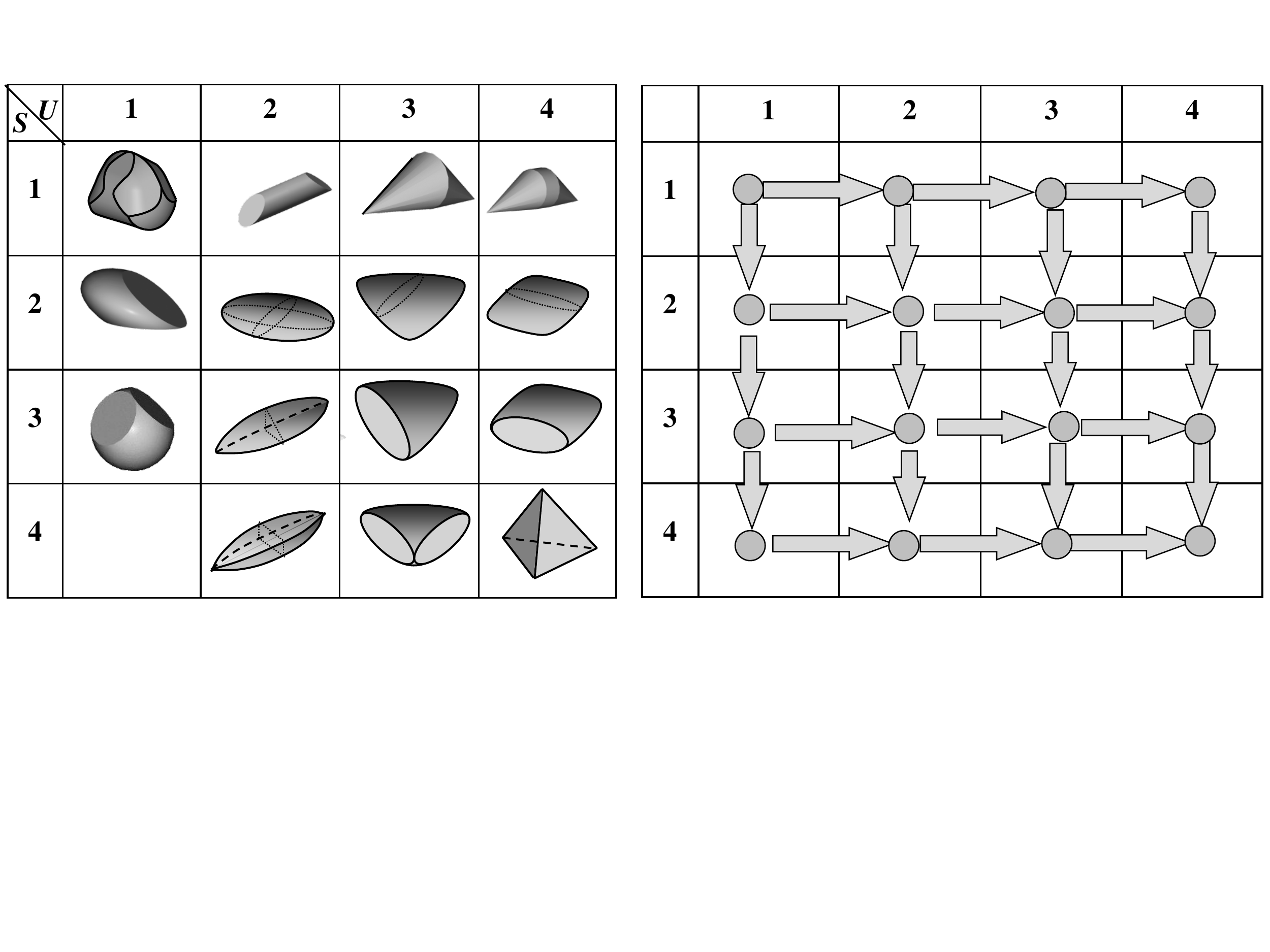}
\vspace{-3cm}
\caption[]{Primary equilibrium classes. Left: examples of convex
bodies; rows and columns correspond to the numbers $S$ and $U$ of
sinks and sources, respectively. Right: the `Columbus Algorithm'
of \cite{VD1} defines a hierarchy among primary classes. Arrows
indicate arbitrarily small truncations of the convex body, 
creating one additional sink or source and a saddle-point.}
\label{fig:primary}
\end{center}
\end{figure}

The classification schemes introduced in \cite{VD1} are based on these
numbers. Specifically, the \emph{primary class} of a generic convex
body $K$ is defined as the pair of integers $\{S,U\}$. In \cite{VD1}
it was shown that no primary class $\{i,j\}$ is empty and a hierarchy
among these classes was defined via the \emph{Columbus algorithm}. 
Using explicit truncations that remove small portions from $K$ by
slicing along convex surfaces, this algorithm generates a pair of
convex bodies $K' \in \{i+1,j\}$ and $K'' \in \{i,j+1\}$, as shown
in Figure \ref{fig:primary}. Thus, starting from the
g\"{o}mb\"{o}c $\{1,1\}$, every row and column can be populated,
implying that the primary classification is complete in this `static'
sense.

\begin{figure}[here]
\begin{center}
\includegraphics[width=100mm]{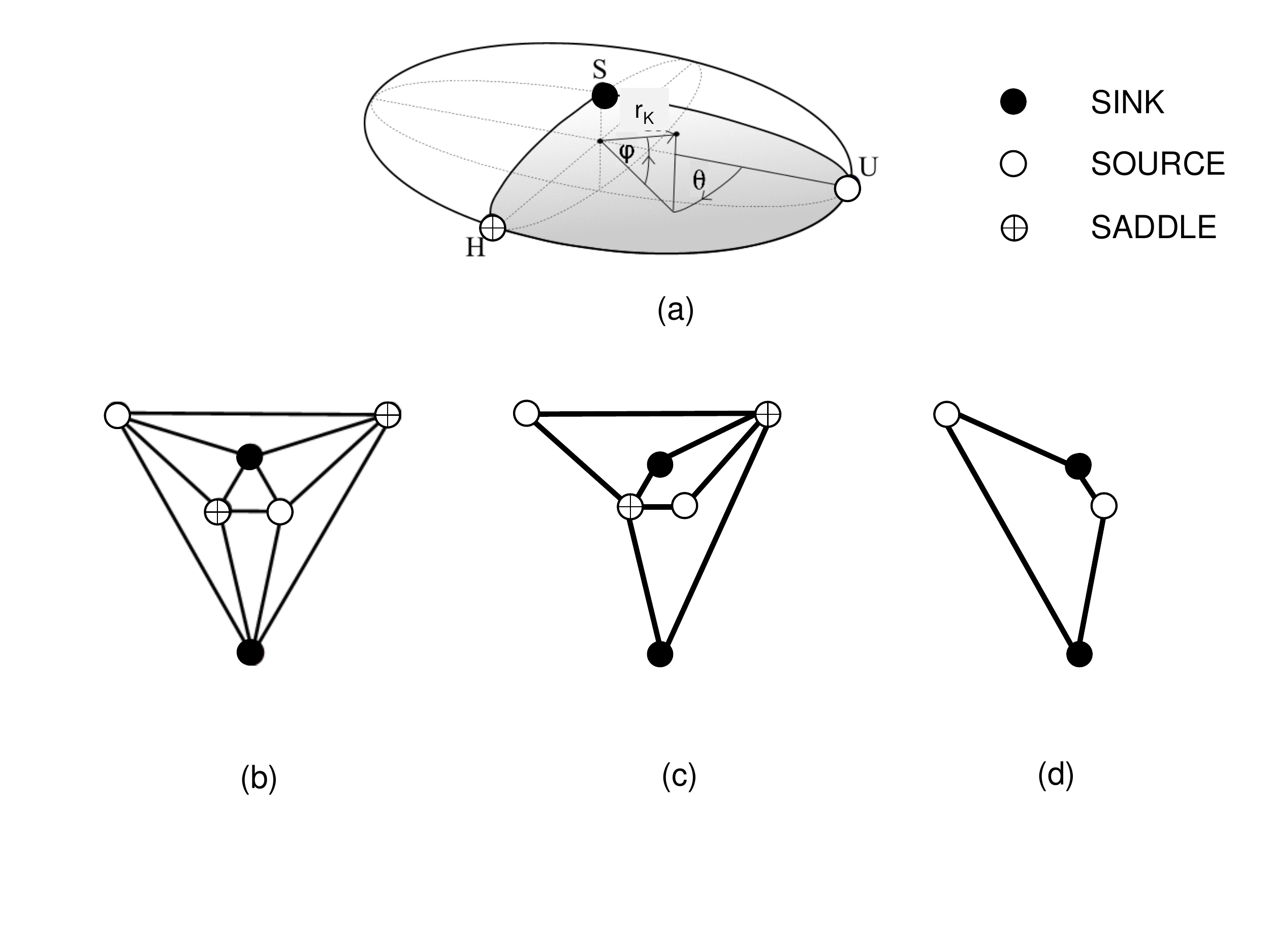}
\vspace{-1cm}
\caption[]{Graph representations of the gradient flow on the tri-axial
ellipsoid in primary equilibrium class $\{2,2\}$. (a) Distance
function $r_K$ given in spherical polar coordinates. (b) 3-colored
quadrangulated  primary representation $Q^3(\mathbf{v})$. (c)
3-colored triangulated representation $T^3(\mathbf{v})$. (d)
Quasi-dual, 2-colored quadrangulated representation
$Q^2(\mathbf{v})$. The colors refer to vertices, identifying them as
sinks, sources, and (in (b,c)) saddles. }
\label{fig:topology}
\end{center}
\end{figure}

More refined methods exist for classifying the properties of gradient
vector fields $\mathbf{v}=\nabla r_K$, including graph representations
of their Morse-Smale complexes \cite{Dong}. The vertices of these
graphs are fixed points of $\mathbf{v}$, and the edges can be either
isolated heteroclinic orbits connecting saddle-points, representative
non-isolated heteroclinic orbits connecting saddles with sinks and
sources, or both. These are called, respectively, the  primary
representation $Q^3(\mathbf{v})$, the triangulated representation
$T^3({\bf v})$, and the quasi-dual representation $Q^2(\mathbf{v})$
Figure \ref{fig:topology} illustrates these representations for the
tri-axial ellipsoid. For brevity, we call all three types the
\emph{topology graphs} associated with $\mathbf{v}$. Note that all
three graphs $Q^3(\mathbf{v}), T^3(\mathbf{v}),$ and $Q^2(\mathbf{v})$
are embedded on $\S^2$ and we will also consider their abstract, 
non-embedded versions $\bar Q^3(\mathbf{v}), \bar T^3({\bf v})$ and
$\bar Q^2(\mathbf{v})$.  We remark that an abstract graph may have
several, non-homeomorphic embeddings in $\S^2$. Precise definitions
will be given in Section \ref{sec:combinatorics}, and these graphs 
will play a key role in the paper.

\begin{figure}[here]
\begin{center}
\includegraphics[width=120mm]{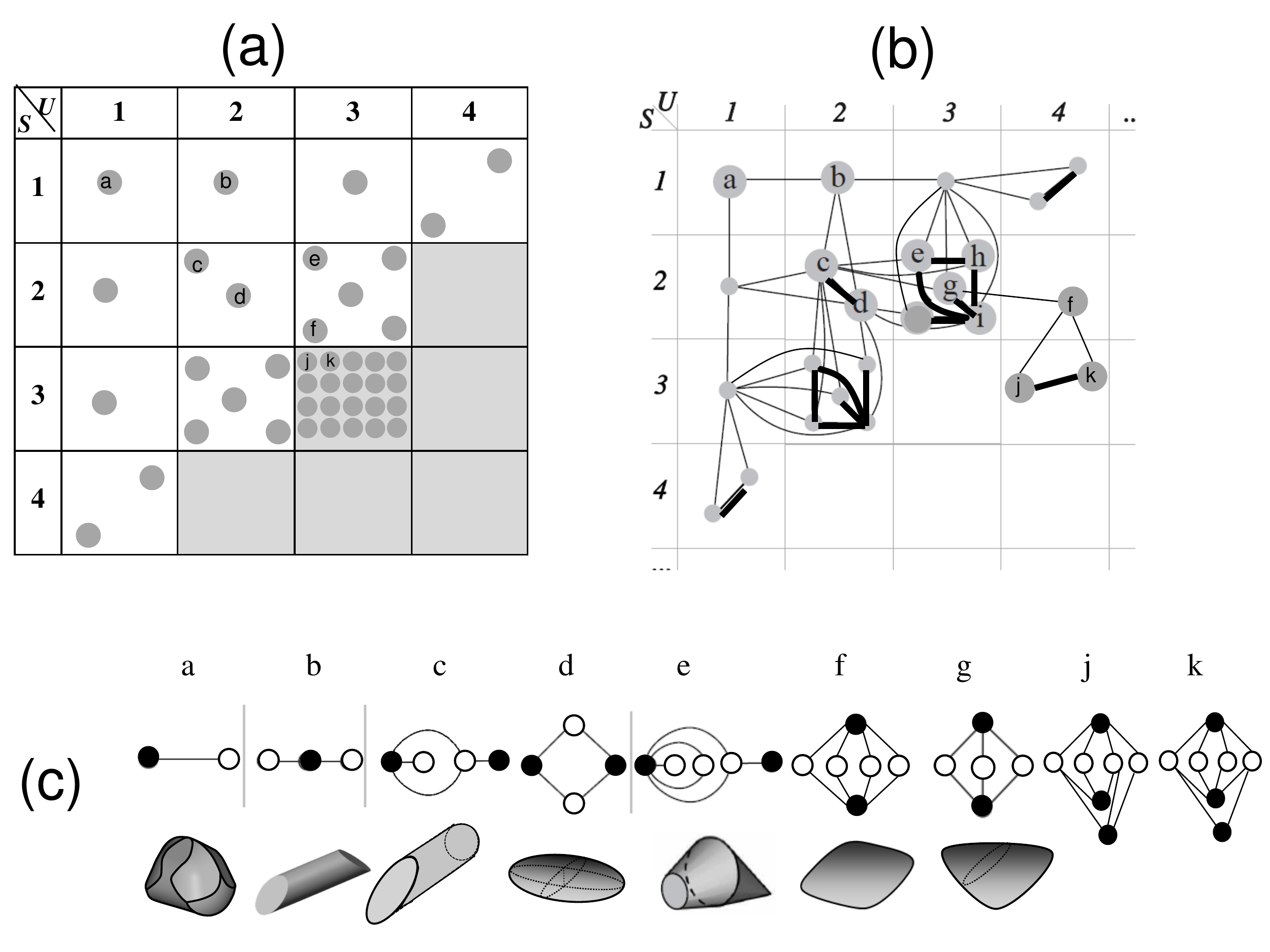}
\caption[]{Secondary and tertiary equilibrium classes. (a) Secondary
and tertiary classes are contained in primary classes. (b) Metagraph
$\mathcal{G}$ with vertices at tertiary classes and edges
corresponding to codimension 1 bifurcations; thin edges: saddle-nodes,
thick edges: saddle-saddle bifurcations. Note that all illustrated
secondary classes contain one tertiary class, so in the figure the
vertices of the metagraph $\mathcal{G}$ correspond simultaneously to
secondary and to tertiary classes. (c)  Quasi-dual topology graphs
$Q^2(\mathbf{v})$ of the tertiary classes labeled a through k in
panels (a,b) and illustrations of a through g as convex bodies }
\label{fig:secondary}
\end{center}
\end{figure}

 We call the class of convex bodies with isomorphic abstract
graphs the  \emph{secondary equilibrium class} and the class of convex
bodies with homeomorphic embedded graphs the \emph{tertiary
equilibrium class} associated with $K$. See Figure
\ref{fig:secondary}(a), which also illustrates that a primary class
can contain different secondary classes: e.g., the ellipsoid is not
alone in class $\{2,2\}$. In \cite{DLSz} it was shown that the
secondary and tertiary schemes are also complete in the sense that no
secondary  or tertiary class is empty.

One can ask whether transitions between different primary, secondary
and tertiary classes are possible within generic families $K(\lambda)$
of smooth convex bodies, parametrized by $\lambda$, as their shapes
change. In generic one-parameter families of gradient vector fields
only two codimension 1 bifurcations occur: saddle nodes and
saddle-saddle connections, and they do so at isolated, critical values
$\lambda = \lambda^{cr}_i$ \cite{GuckenH-83}. Saddle-nodes involve
local changes in topology in which pairs of non-degenerate equilibria,
either a saddle and a sink or a saddle and a source, emerge or
disappear. Saddle-saddle bifurcations are global bifurcations at which
an orbit connecting two saddle-points exists, but the numbers and
types of equilibria do not change. In the former, one of the integers
$S,U$ characterizing the primary class of $K$ increases or decreases
by one; in the latter, the primary class remains unchanged.

\begin{figure}[here]
\begin{center}
\includegraphics[width=120mm]{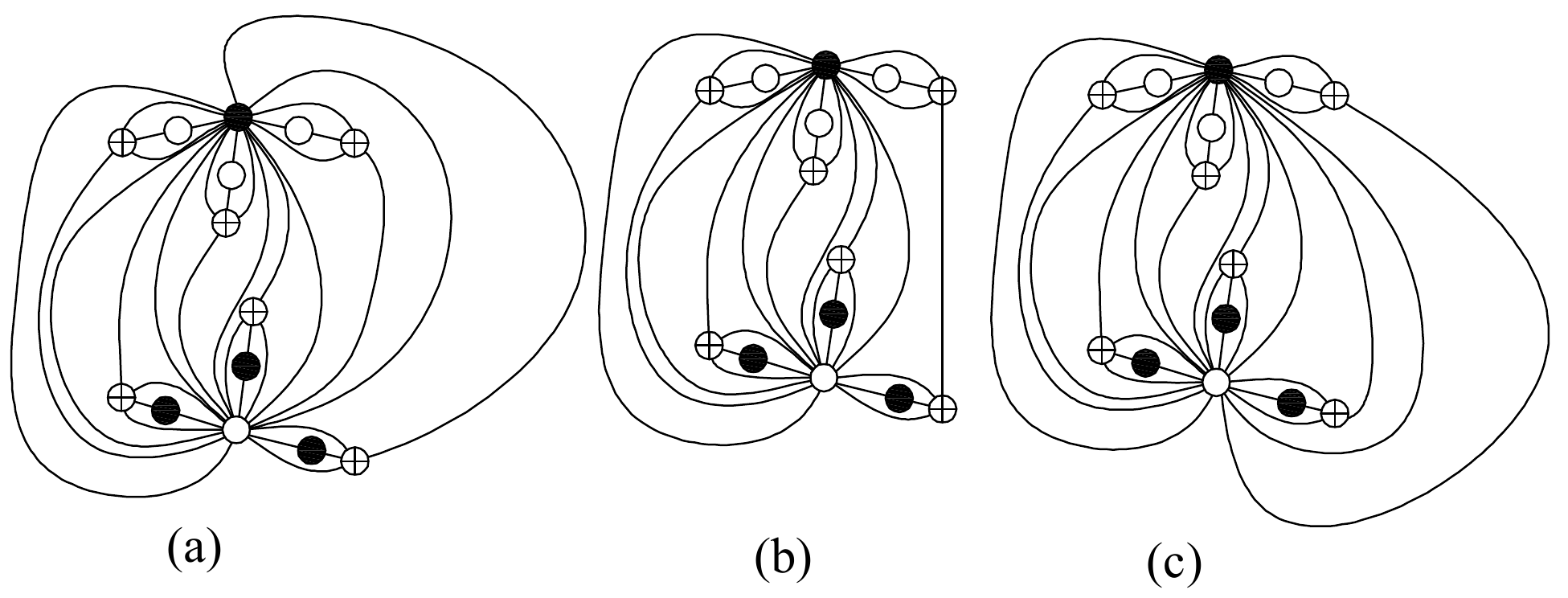}
\caption[]{ Topology graphs (a) and (c) corresponding to two vertices of
$\mathcal G$ in primary equilibrium class $\{4,4\}$ 
and the tertiary edge (b) connecting them. Graphs are shown in
the triangulated representation of graph class $\mathcal{T}^3$. Note
the saddle-saddle connection on (b), and that  (a) and (c)
 are isomorphic as abstract graphs, but not homemorphic as
embedded graphs on $\S^2$.}
\label{fig:tertiary}
\end{center}
\end{figure}

To visualize these transitions we introduce the metagraph
$\mathcal{G}$ with vertices representing the embedded topology graphs
associated with generic gradient fields on $\S^2$. The metagraph
distinguishes \emph{primary edges}, on which saddle-node bifurcations
occur between vertices in different primary classes, from
\emph{secondary edges} that contain saddle-saddle bifurcations between
vertices within the same primary class but in distinct secondary
classes, and \emph{tertiary edges} between vertices in the same
secondary class. In Figure \ref{fig:secondary}(b) primary and
secondary edges are identified by thin and thick lines
respectively. Figure \ref{fig:secondary} shows  only primary
classes with values $S + U \le 6$; here tertiary edges cannot be
illustrated since secondary classes with so few critical points
contain only one tertiary class, and therefore have unique embeddings
on the 2-sphere. Figure \ref{fig:tertiary} illustrates a tertiary edge
 conecting two vertices in the primary class $\{4,4\}$.

Our goal, which we will formally define in subsection
\ref{ss:mainresults}, is to show that not only the vertices, but also
 the primary and secondary edges of $\mathcal{G}$ can be
represented by convex bodies, i.e., physical processes describing
their shape evolution may be represented by paths on $\mathcal{G}$. It
is an intriguing question which of these paths are preferred by
physical abrasion processes.
The investigation of this question lies beyond the scope of the
current paper, but we have some physical intuition on how primary
classifications of convex bodies might evolve. In \cite{D} it is shown
that the evolution of $S$ and $U$ under the partial differential
equations governing collisional abrasion processes can be  modeled by
letting $S$ and $U$ be random variables whose expected values decrease
with time. While this trend has been verified both in laboratory
experiments \cite{Plos} and in the field \cite{Puerto}, almost no
pebbles in the primary classes $\{1,i\},\{j,1\}$, $(i,j=1,2,\dots)$
have been found. In \cite{DL2} a purely geometrical reason for this
phenomenon was pointed out. The difficulty of reducing either the
number of sinks or  sources can be measured by the fraction of volume
that must be removed from a solid, referred to as \emph{robustness}.
In \cite{DL2} it was shown that the robustness of classes
$\{2,i\},\{j,2\}$ is maximal so it is very unlikely that any natural
pebble in these classes will be transformed into any of the classes
$\{1,i\},\{j,1\}$ by natural abrasion.

Much less is known about the secondary, let alone tertiary
classification of convex bodies, and theories for their evolution by
abrasion are lacking. Nevertheless, as for the primary case, field
data indicate that natural shapes strongly prefer some secondary
classes while other classes remain virtually empty. Already in
\cite{Hilbert} it was observed that coastal pebbles tend to be
ellipsoidal. While Rayleigh \cite{Rayleigh1, Rayleigh2,Rayleigh3}
Bloore \cite{Bloore} and Firey \cite{Firey} ultimately showed that the
classical \emph{exact} ellipsoid is not an attracting state in
collisional abrasion processes, nearly ellipsoidal shapes nonetheless
dominate pebble beaches. Without exception, all those shapes in
primary class $\{2,2\}$ for which the secondary classes
were determined have topology graph `d' of Figure 
\ref{fig:secondary}(c), while the other secondary class `c' in
$\{2,2\}$ appears to be missing.  Similar observations apply to other
primary classes. 

As a first step towards understanding these phenomena we show
that the secondary classification scheme of \cite{DLSz} is also
complete in the following `dynamical' sense. Primary and secondary
edges in the metagraph $\mathcal {G}$, containing codimension 1 saddle
node and saddle-saddle bifurcations respectively, exist in the space
of gradient vector fields $\mathbf{v} = \nabla r_K$ on the 2-sphere
associated with convex bodies $K$. In the next subsection we use the
metagraph $\mathcal G$  to formulate our statements and relate them to
earlier results.  Before doing so, we note that the gradient 
vector field $\mathbf{v} = \nabla r_K$ cannot describe the Newtonian
dynamics of the body $K$ rocking on a horizontal plane, which would
require a system of second order differential equations, but that the
stability types of its fixed points correctly reflect those of the 
equilibria of $K$.

\subsection{Definitions and main result}
\label{ss:mainresults}

We first define the metagraph $\mathcal{G}$, whose vertices are
embedded topology graphs representing tertiary equilibrium classes
associated with the Morse-Smale complexes \cite{Milnor} of
gradient vector fields of convex bodies. For simplicity we use the
primary representation $Q^3(\mathbf{v})$, but the triangulated or
quasi-dual representations may also be used to construct
$\mathcal{G}$. The edges of $\mathcal {G}$ correspond to codimension
1 bifurcations connecting these classes, and all possible
one-parameter families of gradient vector fields on the 2-sphere
appear in $\mathcal {G}$. We define the edges and vertices of
$\mathcal{G}$ and we will use these concepts to formulate our results
and relate them to earlier results.

\begin{defn}
Two vector fields $\mathbf{v}$ and ${\bf w}$ on the 2-sphere are
\emph{topologically equivalent} \cite{Sotomayor,GuckenH-83} if
their embedded topology graphs (of the same type) are homeomorphic.
\end{defn}

As noted earlier, in case of generic vector fields \cite{Sotomayor},
topology graphs can be defined by the Morse-Smale complex associated
with the vector field \cite{Dong}. Now we proceed to define the
metagraph $\mathcal{G}$.

\begin{defn}
A vertex of $\mathcal{G}$ is an embedded topology graph
$Q^3(\mathbf{v})$ on $\S^2$, associated with the Morse-Smale
complex of a generic gradient vector field $\mathbf{v}$ on $\S^2$.
\end{defn}

\begin{defn}
The primary class of a vertex $Q^3(\mathbf{v})$ is the pair of
integers $\{S,U\}$, where $S$ and $U$ denote the number of sinks and
saddles of $\mathbf{v}$.  The secondary and tertiary class of a vertex
are the abstract graph $\bar Q^3(\mathbf{v})$  and the embedded graph
$Q^3(\mathbf{v})$, respectively, both associated with the Morse-Smale
complex of $\mathbf{v}$.
\end{defn}

\begin{defn}
An edge of $\mathcal{G}$ is a one-parameter family
$\mathbf{v}(\lambda), \lambda \in [0,1]$ of gradient vector fields
connecting two distinct vertices $Q^3(\mathbf{v}(0))$ and
$Q^3(\mathbf{v}(1))$ of $\mathcal G$. We require that $\mathbf{v}$ is
generic except for  a unique value $\lambda=\lambda^{\star} \in
(0,1)$, for which $\mathbf{v}(\lambda ^{\star})$ exhibits a
codimension 1 bifurcation \cite{GuckenH-83}.
\end{defn}

\begin{defn}
We call an edge $\mathbf{v}(\lambda), \lambda \in [0,1]$ primary if
the primary classes of $Q^3(\mathbf{v}(0))$ and $Q^3(\mathbf{v}(1))$
are different. We call an edge $\mathbf{v}(\lambda), \lambda \in
[0,1]$ secondary if the primary class of $Q^3(\mathbf{v}(0))$ and
$Q^3(\mathbf{v}(1))$ are identical, but their secondary classes are
different. We call an edge $\mathbf{v}(\lambda), \lambda \in [0,1]$
tertiary if both the primary and the secondary class of
$Q^3(\mathbf{v}(0))$ and $Q^3(\mathbf{v}(1))$ are identical.
\end{defn}

\begin{defn}
We call a vertex $Q^3(\mathbf{v})$ of $\mathcal{G}$ physical if
there exists a convex body $K$ such that $\nabla r_K$ is 
topologically equivalent to $\mathbf{v}$.
\end{defn}

\begin{defn}
We call a  primary, secondary or tertiary equilibrium class
physical if it contains at least one physical vertex.
\end{defn}

\begin{defn}
We call an edge $\mathbf{v}(\lambda), \lambda \in [0,1]$ of
$\mathcal{G}$ physical if there exists a one-parameter family
$K(\lambda), \lambda \in [0,1]$ of convex bodies such that 
$\nabla r_K(\lambda)$ is topologically equivalent to
$\mathbf{v}(\lambda)$ for all values of $\lambda \in [0,1]$.
\end{defn}

Now we can formulate earlier and current results. Regarding
primary equilibrium classes we have

\begin{thm} \label{thm1}
All primary classes are physical.
\end{thm}

This result, proved in \cite{VD1}, was generalized in \cite{DLSz}
to include secondary and tertiary classes:

\begin{thm} \label{thm2}
All vertices of $\mathcal G$ are physical.
\end{thm}

In the current paper our goal is to further extend  Theorems 
\ref{thm1} and \ref{thm2} by proving the physical existence 
of an important subset of edges of $\mathcal G$:

\begin{thm} \label{thm:geom}
All primary and secondary edges of $\mathcal G$ are physical.
\end{thm}

\subsection{Sketch of proof}
\label{ss:sketchproof}

As noted above, the local truncations constructed in \cite{DLSz}
modify the Morse-Smale complex of $K$ to produce one-parameter
families of convex bodies in which either $S$ or $U$ is increased by
1. However, these families do not (necessarily) represent edges
in the metagraph $\mathcal G$ since the genericity of the bifurcation
was not guaranteed by the construction in \cite{DLSz}. One-parameter
families connecting vertices at the ends of secondary edges
(saddle-saddle bifurcations) were not even discussed in \cite{DLSz}.
 
Here we extend these results by constructing a 2-parameter family of
convex bodies whose gradient vector fields are generic in the sense
that  certain codimension 1 subsets (curves) in the parameter
plane correspond to vector fields with codimension 1 local saddle-node
and global saddle-saddle bifurcations, forming primary and secondary
edges of $\mathcal G$. We also show that the codimension 1 bifurcation
curves meet in a codimension 2 saddle to saddle-node bifurcation point.

Because secondary edges correspond to codimension 1 global
saddle-saddle bifurcations, the local methods of \cite{DLSz} do not
apply directly. Rather, we achieve our goal in two steps. In
Section~\ref{sec:combinatorics} we prove Combinatorial
Lemma~\ref{lem:comb}, stating that any secondary edge of the metagraph
$\mathcal G$ bounds a triangular face of $\mathcal G$ of which the two
other edges are primary. As shown in Figure~\ref{fig:HH_4} below, the
vertices of such a face represent three topology graphs  that lie in
adjacent primary classes of $\mathcal G$. The triangles (b,c,d) and
(f,j,k)  in Figure \ref{fig:secondary}(b) above provide examples.  We
then appeal to dynamical systems theory \cite{GuckenH-83} in Section
\ref{sec:dynamics} to show that  such a triangular face could
contain a codimension 2 bifurcation point for the gradient flow
$\mathbf{v} = \nabla r_{K(\lambda)}$ and describe how codimension 1
saddle-node and saddle-saddle bifurcations emanate from this point.

In Section~\ref{sec:geometry} we take the second step, providing an
explicit geometrical construction that realizes the codimension 2
bifurcation via an arbitrarily small truncation of $K$ depending on
two parameters. First, in Subsection \ref{ss:one} we prove that a
truncation exists under the assumption that the resulting displacement
of the body's mass center has no effect on the topology of its
gradient flow. Then, in Subsection \ref{ss:two} we construct a
\emph{simultaneous}, auxiliary truncation such that the mass center
remains fixed under the combined truncations, implying that the
topology of the flow is preserved. Finally, in
Section~\ref{sec:summary} we summarize our results and point out some
possible consequences.

\section{Combinatorial part}
\label{sec:combinatorics}

\begin{figure}[ht]
\includegraphics[width=\textwidth]{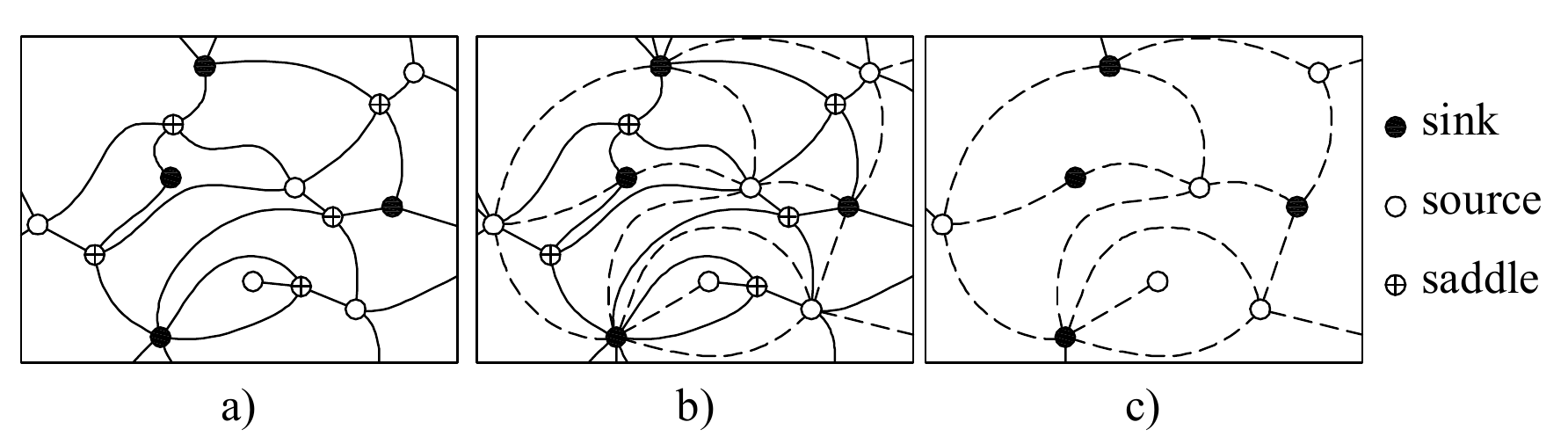}
\caption[]{Different representations of a part of a Morse-Smale
complex. (a) Primary topology graph in class ${\mathcal Q}^{3}$. (b)
Triangulated topology graph in class ${\mathcal T}^{3}$. (c)
Quasi-dual topology graph in class ${\mathcal Q}^2$.}
\label{fig:quasi-dual}
\end{figure}

Before stating the combinatorial lemma, we define three classes of
graphs associated to Morse-Smale complexes on the 2-sphere, of which
the graph representations introduced above and illustrated in
Figure~\ref{fig:topology}(b-d) are examples. As in \cite{DLSz}, we
denote by a \emph{quadrangulation} a finite planar undirected
multigraph on the $2$-sphere in which each face is bounded by a closed
walk of length 4 (cf. \cite{Archdeacon,Brinkmann}). A multigraph
contains no loops but may have multiple (parallel) edges, and it is
usually permitted that the boundary of a face may contain a vertex or
an edge of the graph more than once (e.g. the faces with saddle-source
and source-sink connections in Figure~\ref{fig:quasi-dual}(a,c)). In
addition, we follow Archdeacon et al. \cite{Archdeacon} and regard the
\emph{path graphs} (cf. \cite{Gross}) $P_2$ and $P_3$ as
quadrangulations, where $P_k$ denotes a tree on $k$ vertices, each
with degree at most $2$.

Dong et al. \cite{Dong} introduced three different kinds of graph
to represent a Morse-Smale complex on the 2-sphere, as follows:

\begin{itemize}

\item ${\mathcal Q}^2$ is the class of 2-vertex-colored
quadrangulations. Note that as no quadrangulation contains odd cycles,
each is 2-colorable (cf. \cite{Archdeacon,Nakamoto}). Furthermore, the
coloring of the graph is unique up to switching the colors.

\item ${\mathcal Q}^{3}$ is the class of 3-vertex-colored
quadrangulations with $\deg(p)=4$ for any $p \in H$, and
$|\mathcal{S}|+|\mathcal{U}|-|\mathcal{H}|=2$, where $\mathcal{S}$,
$\mathcal{U}$ and $\mathcal{H}$ denote the sets of vertices of each
given color.

\item ${\mathcal T}^{3}$ is the class of 3-vertex-colored
triangulations with $\deg(p)=4$ for any $p \in \mathcal{H}$, and
$|\mathcal{S}|+|\mathcal{U}|-|\mathcal{H}|=2$, where $\mathcal{S}$,
$\mathcal{U}$ and $\mathcal{H}$ denote the sets of vertices of each
given color.

\end{itemize}

Examples of each class appear in Figure~\ref{fig:quasi-dual}, panels
(c,a,b) respectively.

\begin{figure}[ht]
\includegraphics[width=0.7\textwidth]{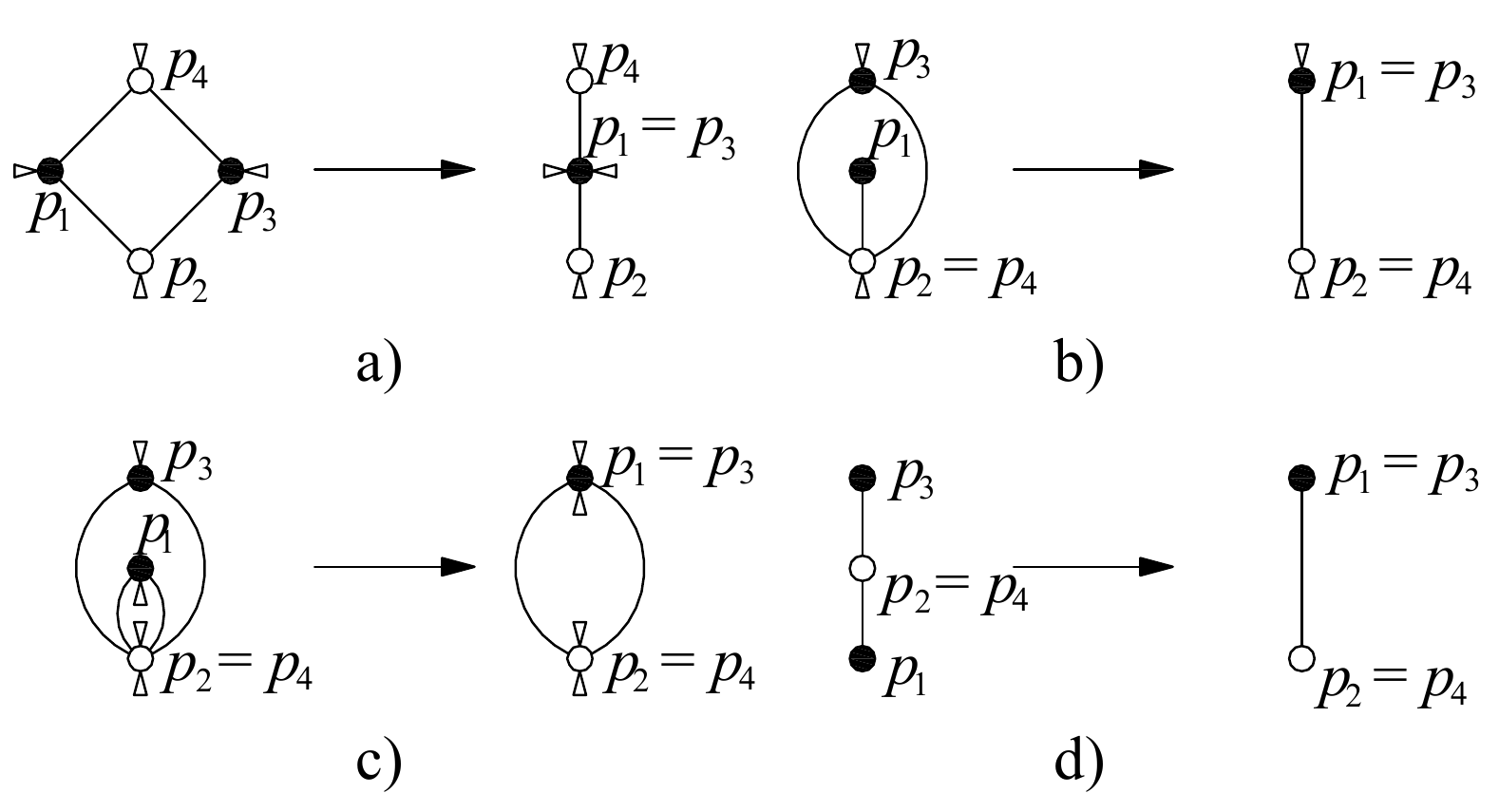}
\caption[]{Face contractions on a graph in class ${\mathcal Q}^2$. As
in \cite{Brinkmann}, triangles incident to some vertices indicate that
one or more edges \emph{may} occur at that position around the
vertex. Here $p_1$ and $p_3$ are sinks; analogous face contractions, 
each removing a source from the graph, can be performed by switching
the colors.}
\label{fig:face-contraction}
\end{figure}

It was shown in \cite{Edelsbrunner} (cf. \cite{Zomorodian}) that a
Morse-Smale complex on the 2-sphere can be uniquely represented by a
3-vertex-colored quadrangulation in ${\mathcal Q}^{3}$, where the
vertex colors represent the 3 types of critical points (maxima, minima
and saddles) and edges correspond to stable and unstable manifolds:
isolated integral curves that end and start at saddle points. Each
quadrangle is bounded by a closed walk consisting of a source, a
saddle, a sink and a saddle in cyclic order around the face and, and
every saddle has degree 4; see Figure~\ref{fig:face-contraction}(a).
Following Dong et al.~\cite{Dong} we call this the \emph{primary}
topology graph.

Saddle points can be removed from the primary graph without losing
information \cite{Dong}: first we connect sources and sinks inside
each quadrangle, producing a \emph{triangulated} topology graph in
class ${\mathcal T}^{3}$; we then remove all saddle points and edges
incident to them, as in Figures~\ref{fig:quasi-dual}(b,c). Since
non-degenerate saddles have degree 4, the resulting graph is a
2-vertex-colored quadrangulation in class ${\mathcal Q}^2$: the
\emph{quasi-dual} topology graph (cf. \cite{Dong}). Here we use the
latter; however, in Section~\ref{sec:geometry}, the primary graph
representation is preferable. All three representations are equivalent
in the sense that they are mutually uniquely identified.

Let $F=(p_1,p_2,p_3,p_4)$ be a face of any $Q \in \mathcal{Q}^2$
(cf. Figure~\ref{fig:face-contraction}(a), left). Pairs of vertices,
and/or edges connecting them, may coincide. Nonetheless, a quasi-dual
representation admits only two kinds of coincidences: two diagonally
opposite vertices, say $p_2$ and $p_4$ may coincide, and in this case
two consecutive edges, say $(p_4,p_1)$ and $(p_1,p_2)$ may coincide:
these two cases are illustrated in
Figures~\ref{fig:face-contraction}(b) and (c), left. Note that in
Figure~\ref{fig:face-contraction}(c) the internal domain bordered by
the edges $(p_4,p_1)$ and $(p_1,p_2)$ is not a quadrangular face and
necessarily contains at least one additional vertex, as indicated by
the triangles in the inner region.
Figure~\ref{fig:face-contraction}(d) shows the remaining two
exceptional cases: the trees $Q = P_3$ and $Q = P_2$. 

The algorithm in \cite{DLSz} is based on repeated application of a
combinatorial graph operation called \emph{face contraction} (cf.
\cite{Archdeacon,Brinkmann} or \cite{Negami}). Applied to the face
$F$ defined in the previous paragraph, this operation results in the
contraction of the vertices $p_1$ and $p_3$ into the same vertex, and
the disappearance of $F$; the modified graphs, depending on the
`shape' of the original face $F$ are shown on the right of each panel
in Figure~\ref{fig:face-contraction}. The inverse operation of face
contractions is called \emph{vertex splitting}. Combinatorially, for
graphs with at least three vertices it can be defined as follows. Let
$p$ be a vertex of the quadrangulation $Q$, with adjacent edges
$E_1,E_2,\ldots, E_k,E_{k+1}=E_1$ in counterclockwise order, and note
that the other endpoints of some of these edges may coincide.
Choose two, not necessarily distinct edges: $E_x$ and $E_y$. Then we
split $p$ into two vertices $p_1$ and $p_3$, and $E_x$ and $E_y$ into
two pairs of edges $E_{x,1}$ and $E_{x,3}$, and $E_{y,1}$ and
$E_{y,3}$, such that $E_{x,1}, E_{x+1}, \ldots, E_{y-1}, E_{y,1}$ are
connected to $p_1$, and $E_{y,3}, E_{y+1}, \ldots, E_{x-1}$ and
$E_{x,3}$ are connected to $p_3$. This operation can be naturally
modified for primary and triangulated representations: in the primary
representation, instead of two edges we choose two (not necessarily
distinct) faces, whereas in a triangulated  representation we
choose two edges connecting a sink and a source.

In a quasi-dual graph, a transition via a saddle-saddle connection can
be realized as a \emph{diagonal slide} \cite{Negami}, defined as
follows: consider two faces $(v_1,v_2,v_3,v_4)$ and
$(v_4,v_3,v_5,v_6)$ of the quadrangulation sharing an edge
$(v_3,v_4)$, and replace this edge either by $(v_1,v_5)$ or
$(v_2,v_6)$. Then the two faces
$(s_1,u_3,s_2,u_2)$ and $(u_1,s_3,u_2,s_2)$ are replaced by
$(s_1,u_3,s_3,u_2)$ and $(u_1,s_3,u_3,s_2)$. 

To formulate the lemma, we need the following definition.

\begin{defn}
\label{defn:twinsplitting}
Let $Q \in {\mathcal Q}^2$ be a quadrangulation. Two vertex
splittings $W$ and $W'$ of $Q$ are called \emph{twin}, if:
\begin{itemize}
\item The same vertex $p$ is split.
\item Let $A_1$ and $A_2$ denote the sets of edges connected to the
two split vertices in $W$, and define $A_1'$ and $A_2'$ similarly for
$W'$.  Then $A_1$ differs in exactly one element from $A_1'$ or
$A_2'$.
\end{itemize}
In this case $Q$ is called the \emph{ancestor} of the two split
graphs. This definition can also be naturally interpreted for
primary and triangulated representations.
\end{defn}

Note that the second property in Definition~\ref{defn:twinsplitting}
implies that $A_2$ also differs in exactly one element from $A_1'$ or
$A_2'$. The graphs in columns $B$ and $C$ of rows (a) and (b) in
Figure~\ref{fig:combinatorics} can be obtained from the graph in
column $A$ of the same row via twin vertex splittings, but the 
graphs in columns $B$ and $C$ of Figure~\ref{fig:combinatorics}(c)
are isomorphic and hence have no ancestor graph. This corresponds to
the degenerate case 3 in the Proof of the  Combinatorial 
Lemma~\ref{lem:comb} below.


\begin{lem}[Combinatorial Lemma]
\label{lem:comb}
Let $B,C  \in {\mathcal Q}^{2}$ be embeddings of two distinct abstract
graphs $\bar B, \bar C$ in $\S^2$, respectively, such that there is a
diagonal slide that transforms $B$ into $C$. Then there is an
embedding  $A \in {\mathcal Q}^2$ and a pair of \emph{twin} vertex
splittings $W_B$ and $W_C$ of $\bar{A}$ such that $W_B$ transforms $A$
into $B$, and $W_C$ transforms $A$ into $C$.
\end{lem}

We remark that, as we will see in the proof of Lemma~\ref{lem:comb},
there are diagonal slides between non-homeomorphic drawings of the
same graph which cannot be derived from the same ancestor via twin
vertex splittings. Note also the essential condition that the
abstract graphs $\bar B, \bar C$ should be non-isomorphic; this
condition excludes tertiary edges from our argument.

\begin{proof}

To simplify the proof, we use the triangulated variants of $B$ and
$C$, which with a little abuse of notation, we also denote by $B$ and
$C$. Let the two saddles that are connected by the saddle-saddle
bifurcation be denoted by $h_1$ and $h_2$. This edge belongs to two
faces of $B$, say $(s_2,u_1,h_1)$ and $(s_2,u_1,h_2)$, and similarly,
two faces of $C$, say $(s_1,u_2,h_1)$ and $(s_1,u_2,h_2)$.  We note
that, due to the degeneracy of the graph, some of the vertices or
edges may coincide; nevertheless, due to the saddle-saddle connection,
$h_1$ and $h_2$ are distinct.

\begin{figure}[here]
\includegraphics[width=0.7\textwidth]{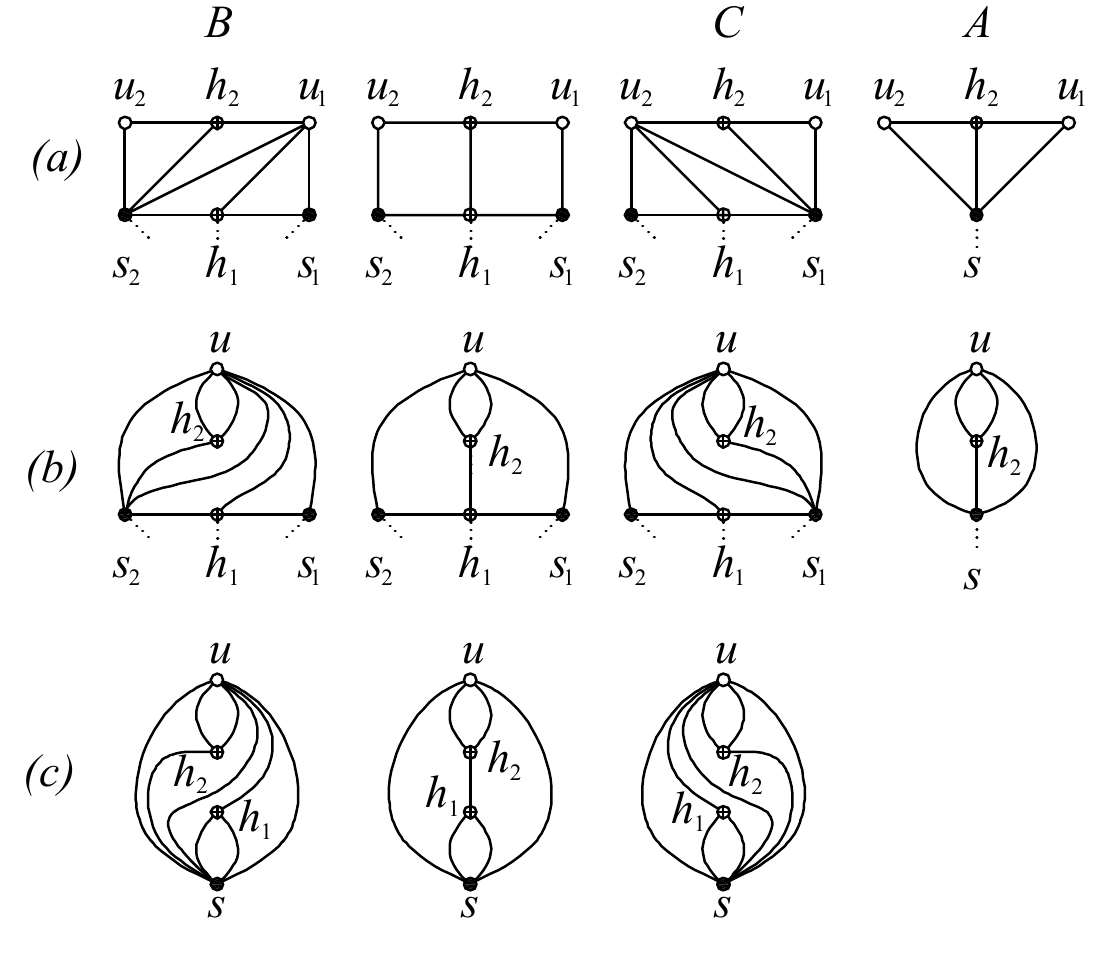}
\caption[]{The connection between diagonal slides and twin vertex
splittings for $\mathcal{T}^{3}$-class graphs. The graphs in row (c),
columns $B$ and $C$ are isomorphic, and thus have no common ancestor,
hence no graph appears in column $A$  of row (c). See the proof
of Lemma~\ref{lem:comb} for further details.}
\label{fig:combinatorics}
\end{figure}

\emph{Case 1}, $s_1$ and $s_2$, and also $u_1$ and $u_2$ are
distinct. Figure~\ref{fig:combinatorics}  row (a) shows the
corresponding faces of $B$, the saddle-saddle bifurcation, $C$, and
the common ancestor $A$  from left to right. Face contractions are
carried out by collapsing the edges $(s_1,h_1)$ and $(s_2,h_1)$
into a single vertex $s$, and the dotted edges starting at $s_1$,
$h_1$ and $s_2$ are contracted into the single dotted edge of
$A$. Furthermore, the edges $(s_2,u_1)$, $(h_1,u_1)$ and $(s_1,u_1)$
are contracted to $(s,u_1)$ in $B$, whereas $(s_2,u_2)$,
$(h_1,u_2)$ and $(s_1,u_2)$ are contracted  to $(s,u_2)$ in
$C$. Since $(s,u_1)$ and $(s,u_2)$ are consecutive edges of $A$ in the
quasi-dual representation, the vertex splittings belonging to the two
face contractions are indeed twin. We remark that in Case 1 another
ancestor can be found by contracting $(u_1,h_2)$ and $(u_2,h_2)$.

\emph{Case 2}, exactly one of the pairs $\{ s_1 , s_2 \}$ or
$\{u_1,u_2 \}$ coincide. Without loss of generality, we may assume
that $u_1 = u_2 = u$ and $s_1 \neq s_2$. Note that as the degree of a
saddle point is $4$, in this case there are two edges starting at
$h_2$ and ending at $u$. Figure~\ref{fig:combinatorics}  row (b)
shows the corresponding faces of $B$, the saddle-saddle bifurcation,
$C$, and the common ancestor $A$  from left to right. Face contraction
is carried out by collapsing the edges $(s_1,h_1)$ and $(s_2,h_1)$
into a single vertex $s$.

\emph{Case 3}, $s_1=s_2=s$, and $u_1=u_2=u$. In this case $B$ and $C$
are isomorphic graphs: Figure~\ref{fig:combinatorics}  row (c). We
note that in this case the two edges starting at $s$ and ending at $u$
may also coincide.

\end{proof}

\section{Dynamical part} 
\label{sec:dynamics}

In this section we describe how codimension 1 saddle-node and
saddle-saddle bifurcations can meet in a codimension 2 bifurcation
of a gradient vector field $\mathbf{v}$ on $\S^2$. Such a bifurcation
point can be associated with each triangular face of the metagraph
$\mathcal G$ having two primary edges and one secondary edge.
We construct an explicit polynomial function $V_{\mu_1,\mu_2}(x,y)$,
depending on two parameters $\mu_1, \mu_2$, that captures the behavior
of $\mathbf{v}$ near a degenerate saddle-node whose strong stable
manifold contains one branch of the unstable manifold of a
non-degenerate (hyperbolic) saddle point. The parameters $\mu_1,
\mu_2$  provide local coordinates on the face of the metagraph near
the codimension 2 point. Since the saddle-saddle or
\emph{heteroclinic} connection is a global phenomenon, our vector
field will necessarily be non-local, but we can nonetheless find a
cubic potential function that captures the local saddle-node and the
global heteroclinic connection. A \emph{homoclinic} orbit to a
saddle-node bifurcation point was previously shown to occur in the
averaged equations for the periodically forced van der Pol oscillator
\cite{HolmesRand-vdP78}, cf. \cite[Sec.~2.1,Figs.~2.1.2-3]{GuckenH-83}.

We first recall the normal form of an isolated codimension 1
saddle-node in a gradient vector field on the plane, which can
be described by a potential function depending on one parameter
\cite{GuckenH-83}:
\ben \label{eq:snV}
V_{\mu_1}(x,y) = \frac{x^3}{3} + \frac{y^2}{2} - \mu_1 x ,
\een
The corresponding vector field
\ben \label{eq:snode}
\brr{l}
\dot{x} = - \frac{\partial V_{\mu_1}}{\partial x} = - x^2 + \mu_1 , \\
\dot{y} = - \frac{\partial V_{\mu_1}}{\partial y} = - y ,
\err
\een
has no fixed points for $\mu_1 < 0$, a saddle-node at $(x,y) = (0,0)$
for $\mu_1 = 0$, and a hyperbolic saddle and a sink at $(x,y) = (-
\sqrt{\mu_1},0)$ and $(+\sqrt{\mu_1},0)$ respectively for $\mu_1 > 0$.

\begin{figure}[here!]
\includegraphics[width=0.85\textwidth]{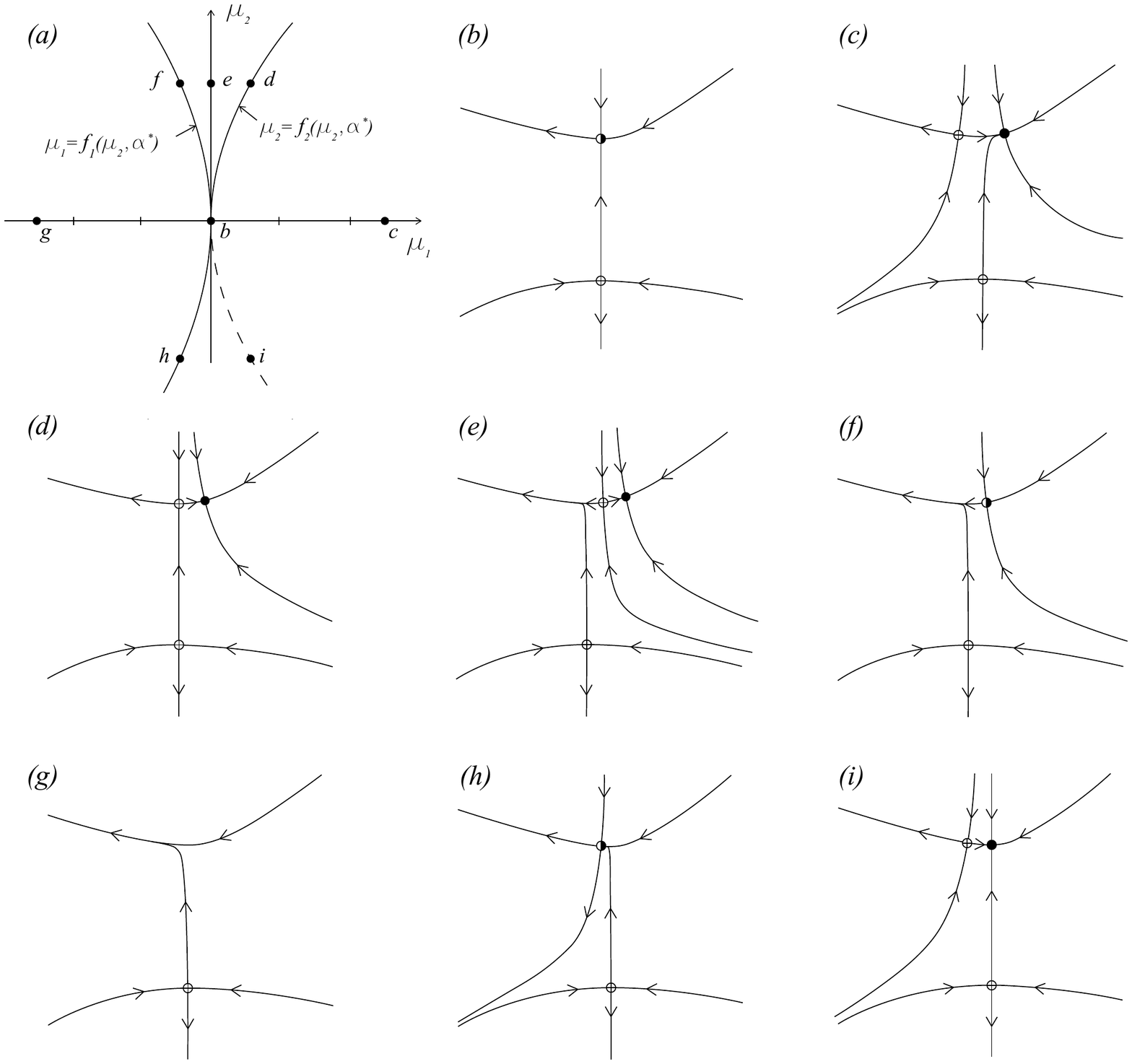}
  \vspace{-4.5cm}
\caption[]{Bifurcations of the gradient vector field
\eqref{eq:cod2}. (a) The bifurcation set in $(\mu_1, \mu_2)$-space near
the codimension 2 point $(0,0)$: saddle-nodes occur on the curve
$\mu_1 = f_1(\mu_2; \alpha) \neq 0$; saddle-saddle connections exist
on the curve $\mu_1 = f_2(\mu_2; \alpha) = \{ \mu_2^2 / 4 \alpha^2 |
\mu_2 \in (0, \sqrt{\alpha}) \}$. (b) At $\mu_1 = \mu_2 = 0$ the
unstable manifold of the saddle at $(0,-1)$ lies in the strong stable
manifold of the saddle-node at $(0,0)$. (c) For $\mu_1 > {\rm max} \{
f_1(\mu_2; \alpha), f_2(\mu_2; \alpha) \}$ two unconnected hyperbolic
saddles coexist with a sink. (d) Along $\mu_1 = f_2(\mu_2; \alpha) <
\alpha^2$ with $\mu_2 > 0$ a codimension 1 saddle-saddle connection
exists on the invariant line $x = -\sqrt{\mu_1}$. (e) For $f_1(\mu_2;
\alpha) < \mu_1 < f_2(\mu_2; \alpha)$ and $\mu_2 > 0$ two hyperbolic
saddles and a sink exist with the lower saddle's unstable manifold
passing left of the upper saddle's stable manifold. (f) On $\mu_1
= f_1(\mu_2;  \alpha) < 0$ with $\mu_2 > 0$ a codimension 1
saddle-node occurs with the lower saddle's unstable passing to its
left. (g) For $\mu_1 <  f_1(\mu_2; \alpha)$ only the lower saddle
exists. (h) On $\mu_1 =  f_1(\mu_2; \alpha) < 0$ with $\mu_2 < 0$ a
codimension 1 saddle-node occurs with the lower saddle's unstable
entering from its right. (i) Along $\mu_1 = f_2(\mu_2; \alpha)$ with
$\mu_2 < 0$ the lower saddle's unstable manifold lies on the invariant
line $x = -\sqrt{\mu_1}$ and intersects the strong stable manifold of
the sink; this is \emph{not a bifurcation point}.}
\label{fig:cod2vf}
\end{figure}

We now add further cubic terms and a linear term containing another
parameter $\mu_2$ to $V_{\mu_1}$ to produce a second hyperbolic
saddle that can be displaced relative to the saddle and sink described
above. We set
\ben \label{eq:SNV}
V_{\mu_1,\mu_2,\alpha}(x,y) = \frac{x^3}{3} + \frac{y^2}{2}
 + \frac{y^3}{3} - \alpha x^2 y - \mu_1 x - \mu_2 x y , \
 \mbox{with} \ \alpha \ge (1/4)^{1/3} ,
\een
so that the vector field \eqref{eq:snode} becomes
\ben \label{eq:cod2}
\brr{l}
\dot{x} = - x^2 + 2 \alpha x y + \mu_1 + \mu_2 y , \\
\dot{y} = \alpha x^2 - y - y^2 + \mu_2 x .
\err
\een
Elementary calculations and linearization at the fixed points show
that, for $\mu_1 = \mu_2 = 0$, the saddle node remains at $(0,0)$ and
a hyperbolic saddle lies at $(0,-1)$. Moreover, the $y$-axis is an
invariant line, because $\dot{x} \equiv 0$ for any solution with
initial condition $(0,y_0)$. The unstable manifold of the saddle
$(0,-1)$ is the line segment $\{ x=0 | y \in (-\infty, 0) \}$, the
upper part of which coincides with the the lower part of the strong
stable manifold  $\{ x=0 | y \in (-1, +\infty) \}$ of the saddle-node.
A disk containing these two fixed points constitutes  a chart,
containing the codimension 2 degenerate vector field, that can be
mapped onto a an open set of $\S^2$: see Figure~\ref{fig:cod2vf}(b).
The term $-\alpha x^2 y$ is necessary to make the lower saddle
hyperbolic (its eigenvalues are $-2 \alpha$ and $+1$). A second
hyperbolic saddle lies  at $(2 \alpha/(4 \alpha^3 - 1), 1/(4 \alpha^3
- 1)$, but this fixed point is irrelevant to the bifurcations of
interest, and it can be driven out of any compact region by letting
$\alpha \rightarrow (1/4)^{1/3} \eqdef \alpha^* \approx
0.62996$. For the cases shown in Figure~\ref{fig:cod2vf}(b-j) we set
$\alpha = 0.62996$.

We now describe the codimension 1 bifurcations and structurally stable
vector fields that emerge from the codimension 2 bifurcation point for
small $\mu_1, \mu_2$. Setting $\dot{x} = \dot{y} = 0$ in \eqref{eq:cod2},
and noting that $y = (x^2 - \mu_1)/(2 \alpha x + \mu_2)$ from the first
equation, we may eliminate $y$ from the second equation to obtain the
fixed point condition
\ben \label{eq:fp1}
F_{\mu_1,\mu_2,\alpha}(x) = a_4 x^4 + a_3 x^3 + a_2 x^2 + a_1 x + a_0 = 0 , 
\een
where
\bea
&& a_4 = 4 \alpha^3 - 1, \ a_3 = 2 \alpha (4 \alpha \mu_2 - 1), \
 a_2 = 2 \mu_1 + 5 \alpha \mu_2^2 -\mu_2,  \nonumber \\
&& a_1 = 2 \alpha \mu_1 + \mu_2^3 \ \ \mbox{and} \ \
 a_0 =\mu_1 \mu_2 - \mu_1^2 .
\label{eq:fp2}
\eea
For $\mu_1 = \mu_2 = 0$ Eqn.~\eqref{eq:fp2} becomes $((4 \alpha^3 -
1)x - 2 \alpha)x^3 = 0$, with a triply-degenerate root at $x=0$ and
the irrelevant root at $x = 2 \alpha/(4 \alpha^3 - 1)$. Setting 
$\alpha = \alpha^*$ so that the latter root lies at $\infty$, the
quartic polynomial becomes a cubic with discriminant
\ben \label{eq:fp2a}
\Delta = 18 a_0 a_1 a_2 a_3 - 4 a_0 a_2^3 +  a_1^2 a_2^2
 - 4 a_1^3 a_3 - 27 a_0^2 a_3^2 .
\een

To obtain an explicit approximation for the saddle-node bifurcation
curve $\mu_1 = f_1(\mu_2; \alpha^*)$, we consider this special case.
Substituting the expressions \eqref{eq:fp2} into \eqref{eq:fp2a} and
setting $\Delta = 0$ yields a polynomial relating $\mu_1$ and $\mu_2$
for which $F_{\mu_1,\mu_2,\alpha^*}(x_0) =
F_{\mu_1,\mu_2,\alpha^*}^\prime(x_0) = 0$ and one of relevant roots
$x_0$ is multiple. Except for $\mu_1 = \mu_2 = 0$, for which $x_0 = 0$
and $F_{\mu_1,\mu_2,\alpha^*}^{\prime \prime}(0) = 0$, this is a
double root, and it corresponds either to a saddle-node bifurcation,
or to the heteroclinic saddle-saddle connection discussed below.
Expanding $\mu_1$ in integer powers of $\mu_2$ and using the fact
that $a_3 = - 2 \alpha + {\mathcal{O}}(\mu_2)$ to determine the
leading terms, we find the following expression for the saddle-node
bifurcation:
\ben \label{eq:fp3}
\mu_1 = f_1(\mu_2; \alpha^*) = -\frac{\mu_2^4}{4} - \frac{3 \cdot 2^{1/3}
 \mu_2^5}{8} - \frac{5 \cdot 2^{2/3} \mu_2^6}{8}  + {\mathcal{O}}(\mu_2^7) .
\een

As in Eqn.~\eqref{eq:snode}, $\mu_1$ primarily controls the
saddle-node bifurcation, but the second parameter $\mu_2$ shifts the
relative $x$ positions of the upper and lower saddles, allowing a
codimension 1 heteroclinic connection to form with $\mu_1 \neq 0$.
Specifically, along the curve
\ben \label{eq:fp4}
\mu_1 = f_2(\mu_2; \alpha) = \frac{\mu_2^2}{4 \alpha^2} , \ \
  \mbox{with} \ \  \mu_2 \in (0, \sqrt{\alpha}) ,
\een
both saddle points lie on the invariant line $x = -\sqrt{\mu_1}$ and a
connecting orbit from the lower to the upper saddle exists (their $y$
coordinates are $\frac{1}{2}\left[-1 \mp \sqrt{1 - \mu_2^2 / \alpha}
\right]$ respectively). This bifurcation curve is shown in
Fig.~\ref{fig:cod2vf}(a) for $\alpha = \alpha^*$, together with the
saddle-node curve $\mu_1 = f_1(\mu_2; \alpha^*)$ (the latter's
curvature is exaggerated for clarity). Note that the
discriminant $\Delta = 0$ for $\mu_1 = \mu_2^2/ 4 \alpha^{*2}$ since
both saddles have the same $x$-coordinate. A similar invariant line $x
= +\sqrt{\mu_1}$ connects the lower saddle to  the strong stable
manifold of the sink for $\mu_1 = \mu_2^2/ 4 \alpha^2$  with $\mu_2 <
0$, but since saddle-sink connections are structurally stable, no
bifurcation occurs here (Figure~\ref{fig:cod2vf}(j)).

In addition to the degenerate codimension 2 vector field at $(\mu_1,
\mu_2) = (0,0)$  shown in panel (b), panels (c-h) show representative
vector fields on the codimension 1 bifurcation curves and
structurally-stable vector fields in the three open regions in the
bifurcation set of panel (a). For the unfolding parameters used here,
the saddle-node and saddle-saddle connection bifurcation curves meet
in a quadratic tangency at $(\mu_1, \mu_2) = (0,0)$. The geometrical
parameters $(d, \phi)$ chosen in the construction that follows produce
bifurcation curves that meet transversely at $(0,0)$, as shown in
Figure~\ref{fig:HH_4}.

\section{Geometrical part}
\label{sec:geometry}

In this section we prove Theorem~\ref{thm:geom}. To do this, it
suffices to create for any primary or secondary edge $E=\{ v_1, v_2\}$
of the metagraph $\mathcal{G}$ a \emph{suitable,} one-parameter family
$K(\lambda)$ of convex bodies, where $\lambda \in
[\lambda_1,\lambda_2]$, with a unique value $\lambda^\star \in
(\lambda_1,\lambda_2)$ such that the graph of $K(\lambda)$ is
homeomorphic to $v_1$ for any $\lambda \in [\lambda_1,
\lambda^\star)$, homeomorphic to $v_2$ for any $\lambda \in
(\lambda^\star, \lambda_2]$, and to the graph of the $1$-codimension
bifurcation defined by $E=\{ v_1, v_2\}$ at
$\lambda=\lambda^\star$. In this case we can choose a
re-parametrization of this family that will satisfy the topological
equivalence condition of the theorem.

We prove the assertion  only for secondary edges of $\mathcal{G}$,
because for primary edges we may apply a simpler version of the same
argument. As noted before, our argument does not apply to tertiary
edges. Secondary edges correspond to non-local bifurcations, so it is
hard to construct by local truncations a suitable one-parameter family
of convex bodies that corresponds to any given secondary edge. To
ensure that local truncations suffice, we rely on
Lemma~\ref{lem:comb}, stating that any secondary edge belongs to a
triangular face of $\mathcal{G}$ of which the two other edges are
primary. Since the latter correspond to local saddle-node
bifurcations, we can use local truncations. We will show that any
face of $\mathcal{G}$ spanned by two primary edges and one secondary
edge can be realized by a \emph{suitable} 2-parameter family
$K(d,\phi)$ of convex bodies
(cf. Definition~\ref{defn:weaklysuitable}). Such a family has
(among others) the property that it collapses to family described
above if we restrict to any of the three edges of $\mathcal{G}$, so
the existence of this suitable 2-parameter family proves the Theorem.

\begin{figure}[here]
\includegraphics[width=\textwidth]{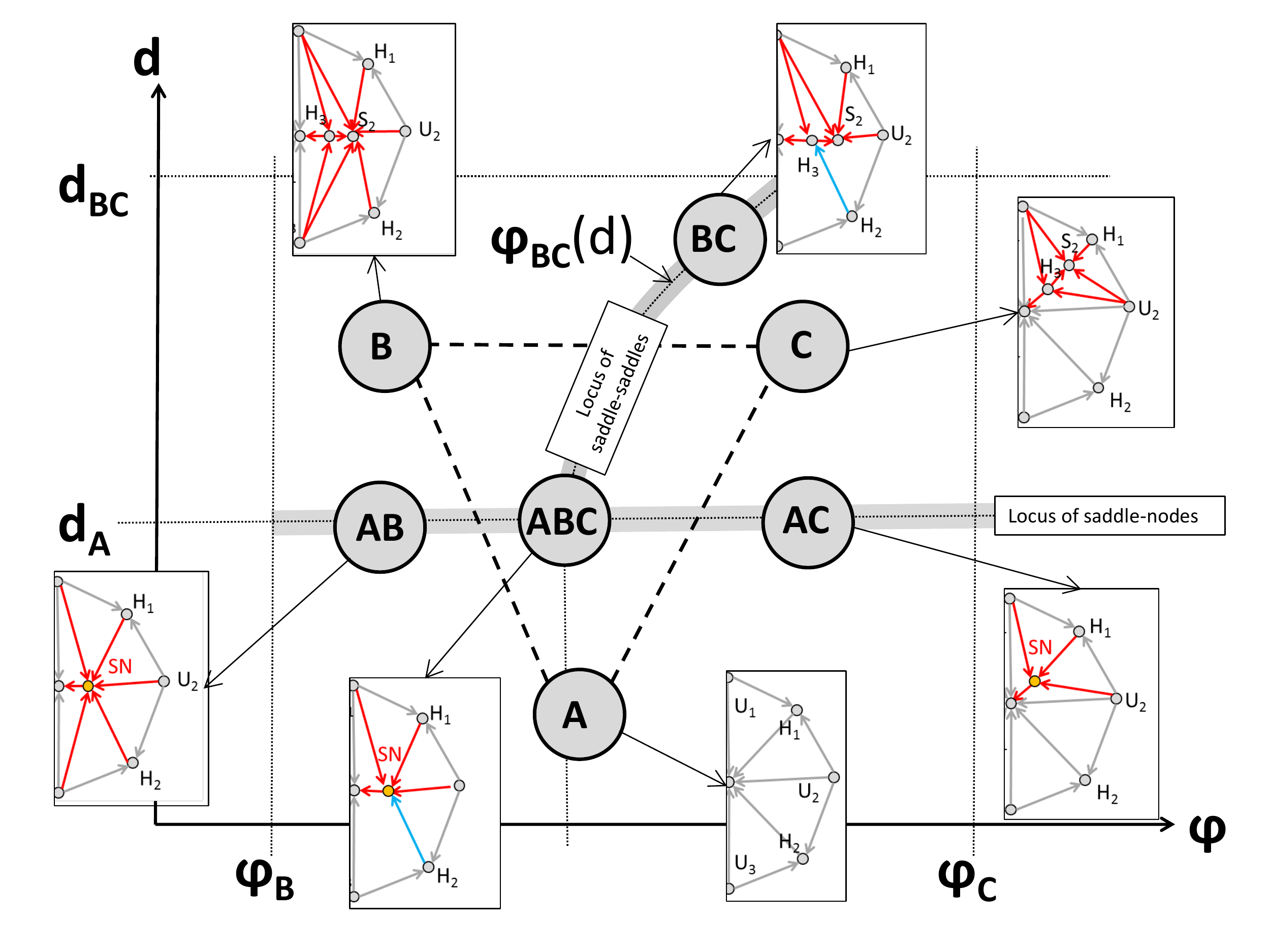}
\caption[]{ A codimension 2 bifurcation on a triangular face of
the metagraph $\mathcal{G}$. The generic graphs $A,B,C$ can be
regarded as subgraphs of a topology graph. For example, by adding one
stable point, they are identical to triangulated representations of
the graphs (f,j,k) in figure \ref{fig:secondary}. The degenerate
graphs $AB$, $AC$, containing codimension 1 saddle-node bifurcations
SN, correspond to primary edges $fj$, $fk$ of $\mathcal{G}$. The 
degenerate graph $BC$, containing a codimension 1 saddle-saddle
connection H2-H3 corresponds to the secondary edge $jk$ of the
metagraph. Finally, the degenerate graph $ABC$, containing the
codimension 2 bifurcation, corresponds to the triangular face $fjk$
of the metagraph $\mathcal{G}$}.
\label{fig:HH_4}
\end{figure}

Let $B$ and $C$ be the primary graph representations of two gradient
vector fields that are connected via any given saddle-saddle
bifurcation. Furthermore, let $A$ be their common ancestor, that is,
$B$ and $C$ can be derived from $A$ by twin vertex splittings. By
Lemma~\ref{lem:comb}, such a graph exists and from \cite{DLSz} we know
that each of the three graphs $A,B,C$ can be associated with the
gradient vector fields of the smooth, convex bodies $K_A,K_B,K_C$,
respectively. We denote the degenerate graphs belonging to the
corresponding transitions by  $AB$, $AC$, $BC$ and $ABC$,
respectively. See Figure~\ref{fig:HH_4}.

\begin{defn}
\label{defn:weaklysuitable}
A 2-parameter family $K(d,\phi)$  of convex bodies, where $d \in
[0,d_{BC}]$ and $\phi \in [\phi _B,\phi_C]$ is called \emph{suitable}
if the function $(d,\phi) \mapsto K(d,\phi)$ is continuous with
respect to Hausdorff distance, and there is a value $d_A \in
(0,d_{BC})$ and a function $\phi _{BC} : [d_A,d_{BC}] \in (\phi
_B,\phi _C)$ such that the following holds:

\begin{enumerate}

\item[(\ref{defn:weaklysuitable}.1)] for every $\phi \in [\phi _B,\phi
_C]$, $K(0,\phi)=K_A$,

\item[(\ref{defn:weaklysuitable}.2)] for every $\phi \in [\phi _B,\phi
_C]$ and $d < d_A$, the graph of $K(d,\phi)$ is homeomorphic to $A$,

\item[(\ref{defn:weaklysuitable}.3)] for every $d > d_A$ and $\phi <
\phi_{BC}(d)$, the graph of $K(d,\phi)$ is homeomorphic to $B$,

\item[(\ref{defn:weaklysuitable}.4)] for every $d > d_A$ and $\phi >
\phi_{BC}(d)$, the graph of $K(d,\phi)$ is homeomorphic to $C$,

\item[(\ref{defn:weaklysuitable}.5)] for every $\phi <
\phi_{BC}(d_A)$, the graph of $K(d_{A},\phi)$ is homeomorphic to $AB$,

\item[(\ref{defn:weaklysuitable}.6)] for every $\phi >
\phi_{BC}(d_A)$, the graph of $K(d_{A},\phi)$ is homeomorphic to $AC$,

\item[(\ref{defn:weaklysuitable}.7)] for every $d > d_{A}$, the graph
of $K(d,\phi _{BC}(d))$ is homeomorphic to $BC$,

\item[(\ref{defn:weaklysuitable}.8)] the graph of
$K(d_{A},\phi_{BC}(d_A))$ is homeomorphic to $ABC$.

\end{enumerate}

If the same properties hold \emph{with the center of mass of $K_A$ as
a fixed reference point}, we say that $K(d,\phi)$ is \emph{weakly
suitable}.
\end{defn}

This definition is illustrated in Figure \ref{fig:HH_4}. We remark
that,  in the context of Section~\ref{sec:dynamics}, the line $\{
d=d_A | \phi \in [\phi _B,\phi_C] \}$ and the curve $\{ \phi =
\phi_{BC}(d) | d \in [d_A,d_{BC}] \}$ form the bifurcation set
associated with the gradient vector field. 

We prove the assertion in two steps: in the first step
(Subsection~\ref{ss:one}), we construct a weakly suitable family.  In
the second step (Subsection~\ref{ss:two}), we modify the construction
in such a way that the center of mass of \emph{every} member of the
family coincides with that of $K_A$, showing that the previously
constructed family is not only weakly suitable but can be made
suitable.

\subsection{Neglecting the motion of the center of mass}
\label{ss:one}

In the first step of the proof we assume that the graph of
\emph{every} convex body is taken with respect to the center of mass
of $K_A$, i.e. we assume that the displacement of the center of mass
does not influence the topology of the flow. For brevity we set
$K=K_A$,  and we consider only the case that the equilibrium point
of $K$ to be split is stable; if it is unstable, a similar argument
with an arbitrarily small, conical extension of the surface can be 
applied.

Let $s$ denote this stable point and the descendant points in the
graphs $B$ and $C$, obtained by splitting $s$, be
$s'_B,s''_B,s'_C,s''_C$,  respectively. Appealing to Lemma 5 of
\cite{DLSz}, we may assume that a  neighborhood of $s$ in $\partial K$
belongs to a sphere $\S$.  Without loss of generality, let the origin
$o$ be the center of this sphere, where the radius of $\S$ is assumed
to be one.  Furthermore, let $c$ denote the center of mass of $K$, and
note that, because $s$ is a stable point, $c$ is contained in the
interior of the segment $[o,s]$.

Let $\Gamma_i$, where $i=0,1,2,\ldots, m$ denote the edges of $A$
starting at $s$, in counterclockwise order around $s$, from outside
$K$. Clearly, for each value of $i$, the part of $\Gamma_i$ in $\S$ is
a great circle arc. These edges are labeled in such a way that the
edges of $B$ starting at $s'_B$ correspond to the $\Gamma_i$'s with
$i=1,2,\ldots,k$ (and those starting at $s''_B$ correspond to the
remaining edges), and the edges starting at $s'_C$ correspond to the
$\Gamma_i$s with $i=1,2,\ldots,k+1$ (and those starting at $s''_C$
correspond to the remaining ones). Observe that, measured in
counterclockwise order, either the angle from $\Gamma_1$ to
$\Gamma_{k+1}$, or the angle from $\Gamma_{k+1}$ to $\Gamma_m$ is less
than $\pi$. Without loss of generality, we may assume that the angle
from $\Gamma_1$ to $\Gamma_{k+1}$ is less than $\pi$.

\begin{figure}[here]
\includegraphics[width=0.4\textwidth]{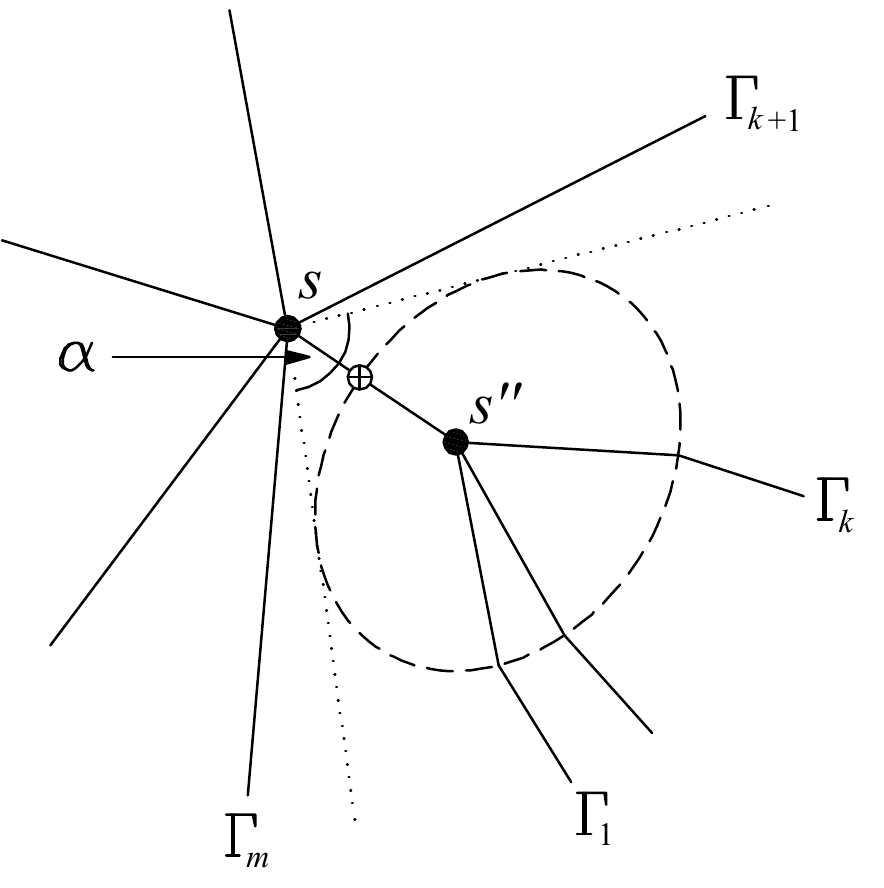}
\caption[]{Truncation by the plane $P(d_{BC},\theta,\phi)$ to create a
body with graph $B$  containing a new saddle on $P \cap \partial
K$ and a sink $s''$.}
\label{fig:geometry1}
\end{figure}

First, we truncate the spherical neighborhood of $s$ by a plane $P$
sufficiently close to but outside $s$, and investigate the equilibrium
points of the truncated body with respect to $c$. In the generic case
we have two possibilities for the graph of the truncated body
$K_P$. If $K \cap P$ does not contain a new stable point, then the
graph of $K_P$ remains homeomorphic to $A$. Furthermore, if $K \cap P$
does contain a new stable point $s''$, then a new saddle point is
created on $P \cap \partial K$, and every heteroclinic orbit on $K$
intersecting $P$ ends up at $s$ (cf. Figure~\ref{fig:geometry1}),
whereas those not intersecting it remain the same.  Finally, note also
that $K \cap P$ contains a stable point if, and only if, the
orthogonal projection of $c$ onto $P$ is contained in the interior of
$K \cap P$. (If the projection is contained on the boundary of this
circle, it is a degenerate case corresponding to a saddle-node
bifurcation.)  We will find a 2-parameter family of planes such
that, if the intersection circle contains the projection of $o$ on any
member, then the edges meeting the circle are either $\Gamma_1,
\ldots, \Gamma_k$, or $\Gamma_1, \ldots, \Gamma_{k+1}$: see Figure
\ref{fig:geometry1}.

An arbitrary plane in $3$-space, and thus, in particular, the
truncating plane, can be defined with three parameters. For this
purpose we use the following coordinates:

\begin{enumerate}

\item[(i)] $d$: the depth of the cut (i.e. the height of the truncated
spherical cap), measured from the point of the sphere where the
tangent plane is parallel to the cutting plane;

\item[(ii)] $\theta$: the arc distance of the center $c_P$ of the
intersection circle (the one created by the cutting plane on the
sphere), from $s$, measured on $\S$;

\item[(iii)] $\phi$: the angle of the great circle arc between $c_P$
and $s$, and from some fixed great circle arc starting at $s$.

\end{enumerate}

Up to a linear transformation, these parameters correspond to the
polar coordinates of the vector pointing from the origin $o$ to its
orthogonal projection onto the truncating plane $P$, where the North
Pole of $\S$ is $s$. Henceforth we denote the plane by $P(d,\theta,\phi)$.

Observe that, measured in counterclockwise order, we have
$\angle(\Gamma_1,\Gamma_k) < \angle(\Gamma_1,\Gamma_{k+1}) <
\angle(\Gamma_0,\Gamma_{k+1}) < \angle(\Gamma_0,\Gamma_{k+2})$. Choose
some angle $0 < \alpha < \pi$ satisfying
\begin{equation}
\label{eq:alpha}
\angle(\Gamma_1,\Gamma_{k+1}) < \alpha < \angle(\Gamma_0,\Gamma_{k+1}).
\end{equation}
Furthermore, for any sufficiently small, fixed value $\theta > 0$,
there is a value $d_{BC}=d_{BC}(\theta,\alpha)$ independent of $\phi$
such that for any plane $P$ with parameters $P(d_{BC},\theta,\phi)$,
$\alpha$ is the angle between the two great circle arcs on the sphere,
starting at $s$ and touching the intersection circle.  Hence, by
(\ref{eq:alpha}), there are some $\phi_B < \phi_{BC} < \phi_C$, with
$\phi_B$ and $\phi_C$ depending only on $\alpha$, and $\phi_{BC} =
\phi_{BC}(d)$ depending on $\alpha$ and $d$, such that

\stepcounter{equation}

\vspace{-1.4cm}

\begin{minipage}[c]{\textwidth}
\parbox[c]{0.03\textwidth}{\vskip2.3cm (\arabic{equation})}
\parbox[t]{0.93\textwidth}{\begin{itemize}

\item for any $\phi \in \left[ \phi_B , \phi_{BC}(d) \right)$ the
plane $P(d_{BC},\theta,\phi)$ intersects $\Gamma_i$ if and only if
$i=1,2,\ldots,k$;

\item the plane $P(d_{BC},\theta,\phi_{BC}(d))$ intersects $\Gamma_i$
if and only if $i=1,2,\ldots,k$, and it is tangent to $\Gamma_{k+1}$;

\item for any $\phi \in \left( \phi_{BC}(d) , \phi_C \right]$ the
plane $P(d_{BC},\theta,\phi)$ intersects $\Gamma_i$ if and only if
$i=1,2,\ldots,k+1$ (cf. Figure~\ref{fig:geometry1}).

\end{itemize} } \hfill
\end{minipage}

Now, consider the one-parameter family $P(d_{BC},\theta,\phi)$, with
$\theta$ fixed and depending only on $\phi \in \left[ \phi_B,
\phi_C\right]$.  If, for any value of $\phi$ in this interval, the
projection of $c$ lies on $K \cap P(d_{BC},\theta,\phi)$, then,
depending on the value of $\phi$, the graph of the body truncated by
the plane is homeomorphic to either $B$ or $C$, or in the degenerate
case to $BC$ (cf. Figure~\ref{fig:geometry1}).  Since we intend to use
local truncations only, we would like to guarantee this property for
any sufficiently small value of $\theta > 0$.  Before proceeding
further, we recall two lemmas from \cite{DLSz}.

\begin{lem}\label{lem:smallradius}
Let $r > |s-c|$ and $\delta > 0$ be arbitrary.
Then there is a convex body $K' \subseteq K$ satisfying the following:

\begin{itemize}

\item[(i)] The graph of $K'$ is homeomorphic to $A$.

\item[(ii)] Denoting the critical point of $K'$ corresponding to $s$
by $s'$, $s'$ has a spherical cap neighborhood in $\partial K'$, of
radius arbitrarily close to $r$.

\item[(iii)] Denoting the integral curve of $K'$ corresponding to
$\Gamma_i$ by $\Gamma'_i$ for every $i$, and by $t_i$ and $t'_i$ the
unit tangent vectors of $\Gamma_i$ and $\Gamma'_i$ at $s'$,
respectively, we have that $|t'_i - t_i | < \delta$.

\end{itemize}
\end{lem}

We note that the same statement is proven in \cite{DLSz} for the case
that $s$ is an unstable point, and the radius of its spherical
neighborhood is arbitrarily close to any given value $0 < r < |s-c|$.

\begin{lem}\label{lem:limit}

Let $C$ be the unit circle in the plane $\Re^2$ with the origin $o$ as its
center, and let $c=(0,\tau)$, where $ \tau > 0$. Let $q_1 =
(\mu_1,\nu_1)$ and $q_2=(\mu_2,\nu_2)$ be two points of $C$ such that
$v_1>0$. 

\begin{itemize}

\item[(i)] If $[q_1,q_2]$ is perpendicular to $[s,q_1]$, then
$\lim \limits_{\mu_1 \to 0} \frac{\mu_2}{\mu_1} = \frac{2 \tau}{1-\tau}$.

\item[(ii)] If the angle of $[q_1,q_2]$ and $[c,q_1]$ is
$\frac{\pi}{2}-C \mu_1$ for some constant $C'$ independent of $\mu_1$,
then $\lim \limits_{\mu_1 \to 0} \frac{\mu_2}{\mu_1} = \frac{2
\tau}{1-\tau} + 2C'$.

\end{itemize}
\end{lem}

\begin{figure}[here]
\includegraphics[width=0.22\textwidth]{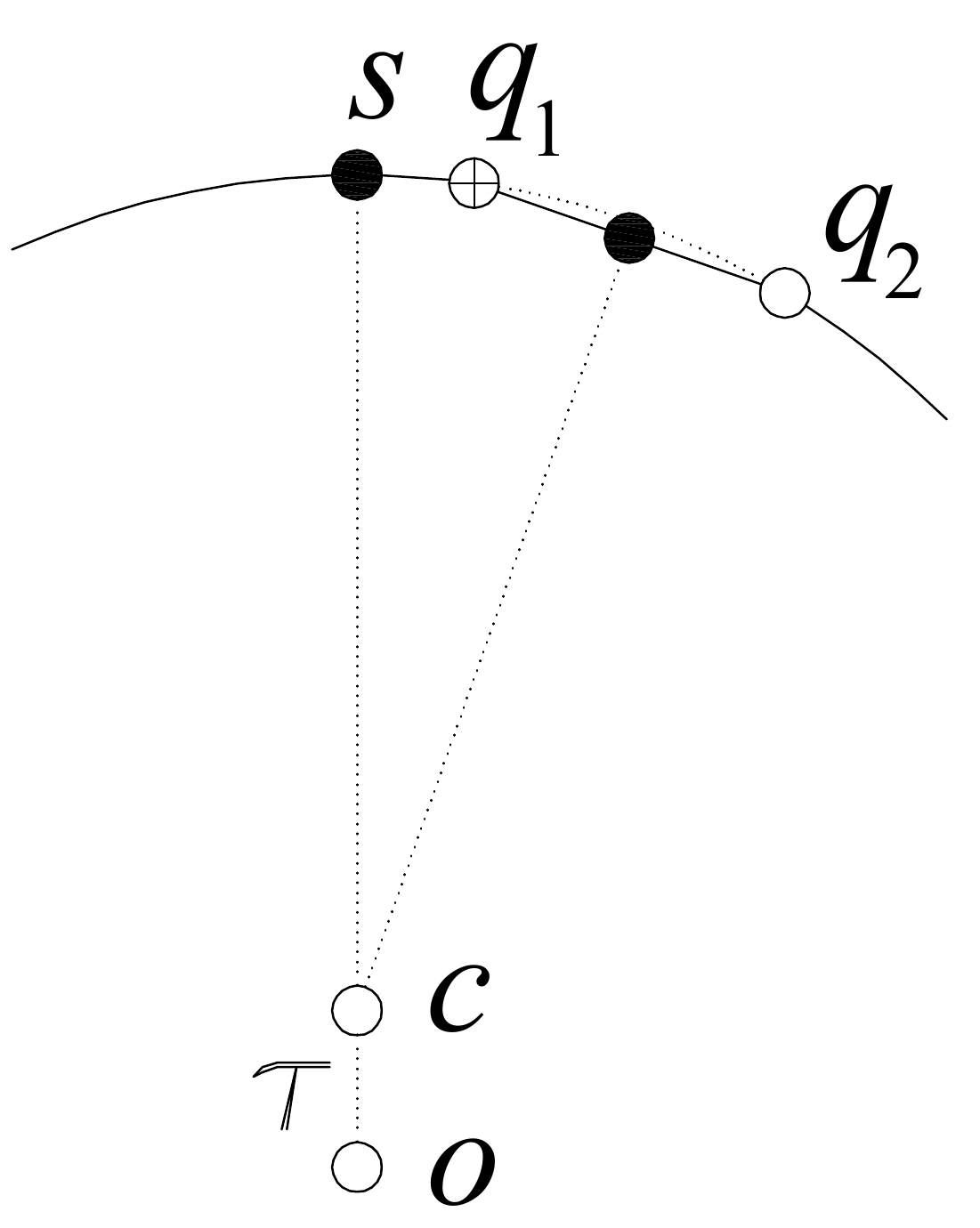}
\caption[]{An illustration for Lemma~\ref{lem:limit}. The sinks
$s$ and $s''$ (filled circles) and the saddle (circle with cross) of
Figure~\ref{fig:geometry1} are identified.}
\label{fig:profileview}
\end{figure}

Now, consider a plane $P=P(d_{BC},\theta,\phi)$ with an arbitrary
value of $\phi$ (cf. Figure~\ref{fig:geometry1}), and let the closest
and the farthest points of the circle $\bar{C} = P \cap \partial K$
from the segment $[o,s]$ be denoted by $q_2$ and $q_1$,
respectively. For convenience, we imagine, for the moment, the plane
containing $q_1 = (\mu_1,\nu_1)$, $q_2 = (\mu_2, \nu_2)$ and
$c=(0,\tau)$ as $\Re^2$ in Lemma~\ref{lem:limit}:
Figure~\ref{fig:profileview}.

For any sufficiently small $\theta > 0$, we require that the
orthogonal projection of $c$ on $P$ lie in the interior of the segment
$[q_1,q_2]$. Since $\angle(q_1,q_2,c) < \frac{\pi}{2}$, this property
holds if and only if $\angle(q_2,q_1,c) < \frac{\pi}{2}$ for any
sufficiently small $\theta > 0$. Recall that $d_{BC}$ is defined by
the fact that the angle of the two tangent lines of $\bar{C}$, passing
through $p$, is equal to $\alpha$.  Let $\bar{C}^s$, $q_1^s$ and
$q_2^s$ denote the central projections of $\bar{C}$, $q_1$ and $q_2$,
respectively, onto the tangent plane of $K$ at $s$.  Then, as $\theta
\to 0$, the limit of the angle of the two tangent lines of
$\bar{C}_s$, passing through $s$, is equal to $\alpha$.  Thus, an
elementary computation yields that, as $\theta \to 0$, the limit of
the ratio of the $x$-coordinate of $q_2^s$ to that of $q_1^s$ is equal
to $\frac{1+\sin \frac{\alpha}{2}}{1-\sin \frac{\alpha}{2}}$, implying
that the same holds for $\lim\limits_{\theta \to 0}
\frac{\mu_2}{\mu_1}$.

 We conclude that our requirement that the orthogonal projection
of $c$ on $P$ lies inside $[q_1,q_2]$ for any sufficiently small
$\theta > 0$ is satisfied if $\frac{2\tau}{1-\tau} < \frac{1+\sin
(\alpha / 2)}{1-\sin (\alpha / 2)}$, but not if $\frac{2\tau}{1-\tau}
 > \frac{1+\sin (\alpha / 2)}{1-\sin (\alpha / 2)}$. To guarantee the
former, we apply Lemma~\ref{lem:limit}, and choose $\delta > 0$
sufficiently small, i.e., such that for the truncated body $K'$ and
heteroclinic orbits $\Gamma_i$, the inequalities (\ref{eq:alpha})
remain true with the same value of $\alpha$. Let $c'$ be the center of
mass of $K'$ and $o'$ be the center of the spherical neighborhood of
$s$. Furthermore, let $\tau'= \frac{|s-c'|}{|s-o'|}$.  Note that for a
suitable choice of $r$,  we have $\frac{2 \tau'}{1-\tau'} >
\frac{1+\sin (\alpha / 2)}{1-\sin (\alpha / 2)}$.

According to the previous paragraph, with a little abuse of notation,
we assume that for the \emph{original body $K$}, for any sufficiently
small $\theta > 0$, the orthogonal projection of $c$ on
$P = P(d_{BC},\theta,\phi)$ lies in the interior of $P \cap K$. Let
$d_A = d_A(\theta,\alpha)$ denote the value of $d$, independent of
$\phi$, at which the projection of $c$ lies on the boundary of $P
\cap K$.

\begin{figure}[here]
\includegraphics[width=\textwidth]{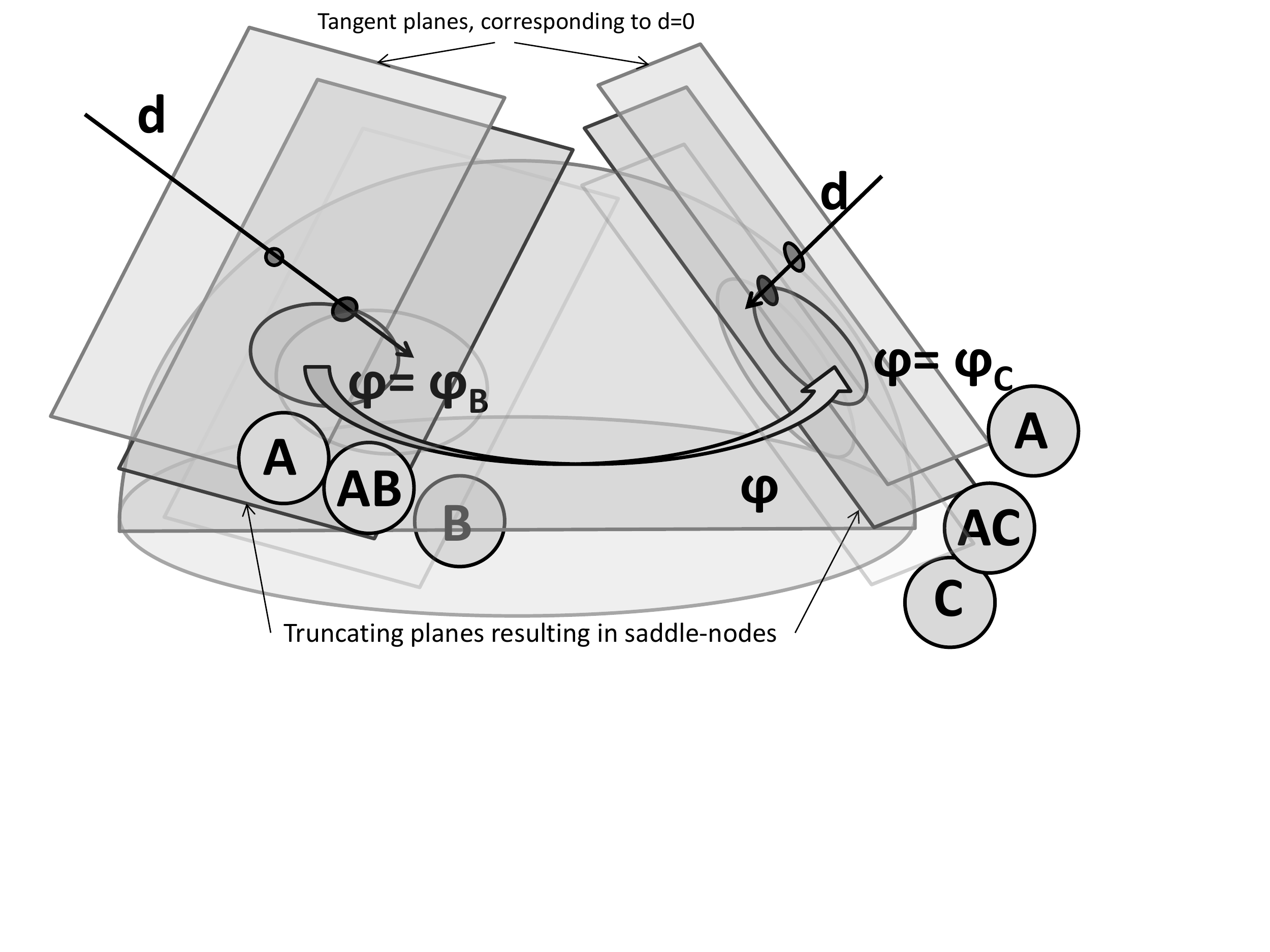}
\vspace{-3cm}
\caption[]{The 2-parameter family of truncations used in the
construction. Variation of the angle $\phi$ of rotation of the
truncating plane results either in a saddle-saddle bifurcation, or in
no bifurcation.  Variation of the depth $d$ of the truncation results
in a saddle-node bifurcation; the graphs belonging to the two extremal
values of $\phi$ are identified with capital letters.}
\label{fig:HH_3}
\end{figure}

We have shown that, for $\phi \in [\phi_B,\phi_C]$ and  $\theta$ is
sufficiently small, the intersection circle $P(d_{BC},\theta,\phi)
\cap K$ contains in its interior a new stable point with respect to
$c$. Thus, the graph of any such truncated body $\bar{K}$ is
homeomorphic to either $B$ or $C$, or to $BC$.  Then we fix a
sufficiently small value of $\theta$, and take the 2-parameter family
of convex bodies $K(d,\phi)$, where $d \in [0,d_{BC}]$, and $\phi \in
[\phi_B,\phi_C]$, defined as the truncation of $K$ by the plane
$P(d,\theta,\phi)$: see Figure~\ref{fig:HH_3}).  Finally, for any
value of $\phi$, $K(0,\phi) = K$, which shows that
(\ref{defn:weaklysuitable}.1) in Definition~\ref{defn:weaklysuitable}
is satisfied. The remaining properties in
Definition~\ref{defn:weaklysuitable} of a weakly suitable family
follow from (3).  This completes the first step of the proof.

We note that the bifurcation diagram of Figure~\ref{fig:HH_4} in
the geometric parameters $d,\phi$ used for the construction of the
truncating plane and that of Figure~\ref{fig:cod2vf}(a) in the
unfolding parameters $\mu_1,\mu_2$ of Section~\ref{sec:dynamics} are
topologically \emph{but not differentiably} equivalent. As noted in
Section~\ref{sec:dynamics}, the bifurcation curves meet in a tangency
in Figure~\ref{fig:cod2vf}(a); however, they meet at a nonzero angle
in Figure~\ref{fig:HH_4}.

\subsection{Annihilating the motion of center of mass by an
auxilary truncation} 
\label{ss:two}

In this subsection we  modify the family $K(d,\phi)$ in such a way
that the center of mass of every member in the modified family 
remains at $c$. To do this we need some additional assumptions on $K$.

Let $L$ be the line passing through $s$ and $c$, and let $w$ denote
the point of $L \cap \partial K$ different from $s$. We show that
$K=K_A$ can be chosen in such a way that $q$ is not an equilibrium
point, and that it does not belong to any edge of $A$. First we modify
the convex body $K_0$ in class $(1,1)$ in \cite{VD1} to satisfy this
property. Since the graph of $K_0$ does not contain edges, we need
only   show that no line through the center of mass passes through
more than one equilibrium point.

Since $K_0$ has $D_4$ rotational symmetry, in a suitable coordinate
system, its two equilibrium points and center of mass $c$ lie on the
$z$-axis, and $K_0$ is symmetric with respect to the
$(x,y)$-coordinate plane. Thus, all the tangent planes of $K_0$,
parallel to the $x$-axis (i.e. satisfying the property that one of
their translates contains the $x$-axis), touch $K_0$ at points in the
$(y,z)$-plane.  Clearly, cutting off a sufficiently small part of
$K_0$ near the positive half of the $x$-axis does not change the
number of equilibria nor the primary equilibrium class  $\{1,1\}$ of
the body. The center of mass $c'$ of the modified body $K_0'$ is in
the open half space $\{ x < 0 \}$.  Hence, if the tangent plane of
$K_0'$ at some point $p$ is perpendicular to the segment $[c',p]$,
then the outer normal vectors of this plane have positive
$x$-coordinates, implying that $p$ is in the open half space $\{ x > 0
\}$. To show that any graph $A$ can be associated to a convex body
$K_A$ satisfying this property, we observe that, by
\cite[Theorem~1]{DLSz}, $K_A$ can be obtained from $K_0'$ by a finite
sequence of local deformations.

\begin{figure}[here]
\includegraphics[width=0.6\textwidth]{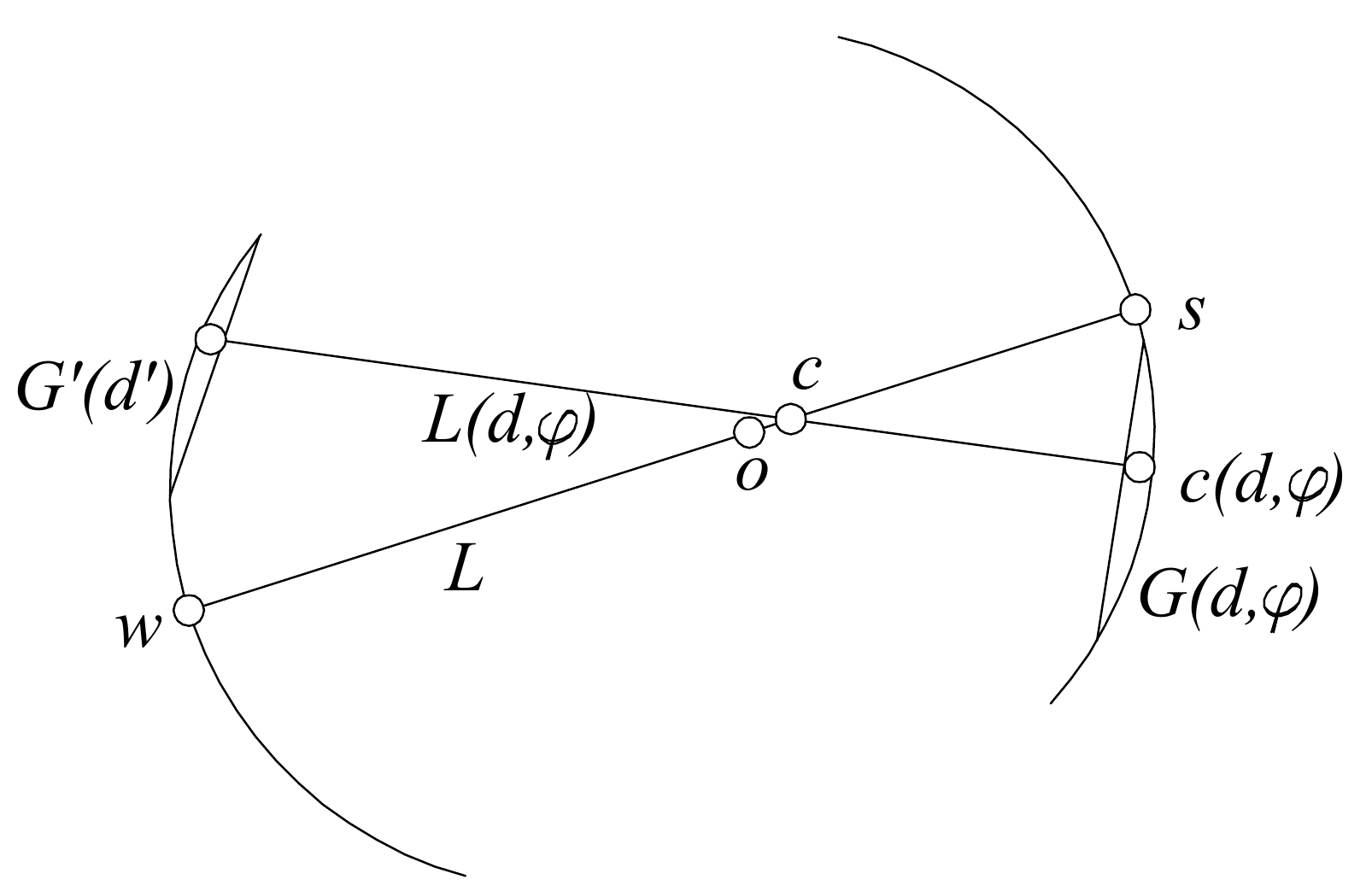}
\caption{The second truncation near the critical point $w$ opposite
to $s$; the circular arcs lie on the spherical caps $G$ and $G'$.}
\label{fig:second_truncation}
\end{figure}

In \cite{DLSz} we also showed that a neighborhood of any point of a
non-isolated heteroclinic orbit, or a sink, or a source can be
truncated by a sphere without changing the class of the graph of the
body.  Furthermore, by \cite[Lemma~1]{DLSz}, we obtain that, applying
a sufficiently small truncation at $w$, the line connecting the
modified stable point and the modified center of mass intersects this
spherical surface.  Thus, we may also assume that a neighborhood of
$w$ is a sphere $\S'$. Nevertheless, note that the center of $\S'$ is
not necessarily on the line $L$. Let $c_{d,\phi}$ denote the center of
mass of the \emph{truncated} spherical cap $G(d,\phi)$ near $s$. To
obtain a modified body $K'(d,\phi)$, we truncate $K_A$ near $q$ by a
second plane $P'(d',\theta',\phi')$, such that the center of mass of
the union of $B$, and the second truncated (open) spherical cap $G'$,
is $c$ (cf. Figure~\ref{fig:second_truncation}). Clearly, in this case
the center of mass of the doubly truncated body $K'(d,\phi)$  is
identical to the center of mass $c$ of $K$.

Let $L(d,\phi)$ denote the line connecting $c$ and $c(d,\phi)$. First,
let $d'$ be fixed. Then, changing $\theta'$ and $\phi'$, the locus of
the centers of mass of $G'$ is a part of a sphere $\S'_{d'}$,
concentric to $\S'$, and the radius of this sphere depends on $d'$ and
$\S'$ only. Thus, if $\theta > 0$ is sufficiently small, for every
line $L(d,\phi)$ and every (small) value of $d'$ there is a unique
position of $G'$ such that its center of mass lies on $L(d,\phi)$. Let
us call this cap $G'=G'(d')$. Note that the center of mass of the
union of $G(d,\phi)$ and $G'(d')$ is $c$ if and only if, the torques
about $c$ exerted by the two caps are equal.  Here, the distance
of the center of mass of $G'(d')$ from $c$ is approximately $|q-c|$;
that is a fixed value. Thus, by continuity, for every pair of values
$d, \phi$, there is at least one value of $d'$ such that the center of
mass of $G(d,\phi) \cup G'(d')$ is $c$.  Let $G'(d,\phi)$ be the
spherical cap $G'(d')$, where $d'$ is the smallest value for which
this property holds. Then, clearly, $G'(d,\phi)$ depends continuously
on $d$ and $\phi$, and the 2-parameter family $K \setminus \left(
G(d,\phi)  \cup G'(d,\phi) \right)$ has the required properties.

\section{Summary} 
\label{sec:summary}

In this paper we showed that the secondary classification of smooth
convex solids, based on the Morse-Smale complexes of their
gradient vector fields, is not only complete in the sense that all
combinatorially possible Morse-Smale complexes can be realized on
smooth, convex bodies, but it is also complete in the more general,
`dynamical' sense that all generic transitions between Morse-Smale
complexes represented by non-isomorphic abstract graphs can be
realized on one-parameter families of convex bodies. Among
trajectories of physical convex shape evolution processes we find
examples of such transitions, so our result implies that from a purely
geometrical viewpoint, there is no restriction on these trajectories.

Theorem \ref{thm:geom} admits only one-parameter families exhibiting
one single bifurcation. However, if we only admit saddle-node
bifurcations then, based on our argument in Section~\ref{sec:geometry}
we can formulate a more general claim.  A codimension one, generic
saddle-node is either a \emph{creation} or an \emph{annihilation},
depending on whether the number of generic critical points
increases or decreases by two.  As stated before, at saddle-saddle
bifurcations the number of generic critical points does not change.

To formulate the claim we introduce

\begin{defn}
A generic, one-parameter family $v(\lambda)$ of gradient vector fields
on the 2-sphere is called  strictly monotone if it contains either
only creations or only annihilations and it does not contain any
saddle-saddle bifurcations.
\end{defn}

Using this concept, we can state the following corollary to 
Theorem~\ref{thm:geom}:

\begin{cor}
\label{cor:geom}
For any generic, strictly monotone, one-parameter family $v(\lambda)$
of gradient vector fields on the 2-sphere there exists a one-parameter
family $K(\lambda)$ of (not necessarily smooth) convex bodies such
that $\nabla r_{K(\lambda)}$ is topologically equivalent to
$v(\lambda)$ for every value of $\lambda$.
\end{cor}

To extend this statement further, we make
\begin{conj}
\label{c1}
Every equivalence class on the family of convex bodies, defined by the
tertiary classification system, is connected. That is, for any two
convex bodies $K_1$ and $K_2$ with the same topology graph $A$ there
is a one-parameter family $K(\lambda)$ of convex bodies, where
$\lambda \in [0,1]$, such that $K(0)=K_1$, $K(1)=K_2$, and the graph
of $K(\lambda)$ is $A$ for every value of $\lambda$.
\end{conj}

If Conjecture \ref{c1} is true, Corollary \ref{cor:geom} can be
extended to include not only strictly monotone, but also generic
families. Although our techniques do not admit the investigation of
tertiary edges of $\mathcal{G}$, we also formulate

\begin{conj}
All tertiary edges of $\mathcal G$ are physical.
\end{conj}

Regarding geophysical applications, we remark that in primary class
$\{2,2\}$ one of the secondary classes  (that of ellipsoids) appears
to be dominant and the other appears to be entirely missing among
natural pebble shapes. Our results show that one \emph{could}
continuously transform members of one class into members of the other
class. Apparently, this process exists in natural abrasion only in one
direction.

\section{Acknowledgments}

This work was supported by OTKA grant T104601 and the J\'anos Bolyai
Research Scholarship of the Hungarian Academy of Sciences. Comments
from Timea Szab\'o are gratefully acknowledged.

\bibliographystyle{abbrv}
\bibliography{References}

\end{document}